\documentclass[12pt]{amsart}
\usepackage{amssymb}

\usepackage{amsbsy}

\newcommand{\vsp}{\vspace{5mm}}
\newcommand{\bt}{\begin{Theorem}}
\newcommand{\et}{\end{Theorem}}

\newcommand{\ei}{\end{itemize}}
\newcommand{\bea}{\begin{eqnarray}}
\newcommand{\eea}{\end{eqnarray}}
\newtheorem{Theorem}{\sc Theorem}
\newtheorem{Lemma}[Theorem]{\sc Lemma}
\newtheorem{Proposition}[Theorem]{\sc Proposition}
\newtheorem{Corollary}[Theorem]{\sc Corollary}
\newtheorem{Definition}[Theorem]{\sc Definition}
\newtheorem{Example}[Theorem]{\sc Example}
\newtheorem{Remark}[Theorem]{\sc Remark}

\newcommand{\be}{\begin{equation}}
\newcommand{\ee}{\end{equation}}

\def\qed{\hfill$\Box$}

\def\CR{{\mathcal {R}}}
\def\CJ{{\mathcal {J}}}
\def\CN{{\mathcal {N}}}
\def\CM{{\mathcal {M}}}
\def\CB{{\mathcal {B}}}

\def\CL{{\mathcal {L}}}
\def\CH{{\mathcal {H}}}
\def\CK{{\mathcal {K}}}

\def\CE{{\mathcal {E}}}
\def\CI{{\mathcal {I}}}
\def\CS{{\mathcal {S}}}

\def\CX{{\mathcal {X}}}

\def\CA{{\mathcal {A}}}

\def\CD{{\mathcal {D}}}
\def\CO{{\mathcal {O}}}
\def\CZ{{\rm\ke}rn.26em
\newcommand\la{{\langle}}
\newcommand\ra{{\rangle}}
\newcommand\lar{\leftarrow}
\newcommand\Lar{\Leftarrow}
\newcommand\rar{\rightarrow}
\newcommand\Rar{\Rightarrow}

\vrule width.02em height.5ex depth0ex \kern.04em \vrule width
.02em height1.47ex depth-1ex \kern-.34em Z}

\def\C{{\rm \kern.24em
 \vrule width.02em
    height1.4ex depth-.05ex
 \kern-.26em C}}

\def\ra{{\rightarrow}}

\def\ei{{\bf e_i}}

\def\oz{{\overline{z}}}
\def\us{{\underline{s}}}
\def\ud{{\underline{d}}}

\def\uQ{{\underline{Q}}}
\def\uX{{\underline{X}}}
\def\uT{{\underline{T}}}
\def\uR{{\underline{R}}}
\def\uV{{\underline{V}}}
\def\uL{{\underline{L}}}
\def\uW{{\underline{W}}}
\def\uS{{\underline{S}}}

\def\uC{{\underline{C}}}

\def\uz{{\underline{z}}}
\def\uD{{\underline{D}}}
\def\uE{{\underline{E}}}

\def\N{{\rm I\kern-.23em N}}
\def\B{{\rm I\kern-.25em B}}
\def\D{{\rm I\kern-.25em D}}
\def\E{{\rm I\kern-.25em E}}
\def\F{{\rm I\kern-.25em F}}
\def\H{{\mathcal H}}
\def\K{{\mathcal K}}
\def\I{{\rm I\kern-.25em I}}

\def\M{{\rm I\kern-.23em M}}
\def\P{{\rm I\kern-.25em P}}
\def\A{{\rm \kern.22em
 \vrule width.02em
    height0.5ex depth 0ex
 \kern-.24em A}}
\def\G{{\rm \kern.24em
 \vrule width.02em
    height1.4ex depth-.05ex
 \kern-.26em G}}
\def\J{{\rm \kern.19em
 \vrule width.02em
    height1.47ex depth 0ex
 \kern-.21em J}}
\def\O{{\rm \kern.24em
 \vrule width.02em
    height1.4ex depth-0.5ex
 \kern-.26em O}}
\def\Q{{\rm \kern.24em
 \vrule width.02em
    height1.4ex depth-.05ex
 \kern-.26em Q}}
\def\S{{\rm \kern.18em
 \vrule width.02em
    height1.4ex depth-.9ex
  \kern.12em
  \vrule width.02em
     height0.7ex depth 0ex
  \kern-.34em S}}
\def\T{{\rm \kern.45em
 \vrule width.02em
    height1.47ex depth 0ex
 \kern-.47em T}}
\def\U{{\rm \kern.30em
 \vrule width.02em
    height1.47ex depth-.05ex
 \kern-.32em U}}
\def\V{{\rm \kern.27em
 \vrule width.02em
    height1.47ex depth-.8ex
 \kern-.29em V}}
\def\W{{\rm \kern.25em
 \vrule width.02em
    height1.47ex depth-0.9ex
 \kern.34em
 \vrule width.02em
    height1.47ex depth-.9ex
  \kern-.63em W}}
\def\X{{\rm \kern.30em
 \vrule width.02em
    height1.4ex depth-1ex
  \kern.12em
  \vrule width.02em
     height0.4ex depth 0ex
  \kern-.46em X}}
\def\Y{{\rm \kern.25em
 \vrule width.02em
    height1.0ex depth 0ex
 \kern-.27em Y}}
\def\Z{{\rm \kern.26em
 \vrule width.02em
    height0.5ex depth 0ex
  \kern.04em
  \vrule width.02em
     height1.47ex depth-1ex
  \kern-.34em Z}}

\begin{document}

\begin{center} {\bf {\Large Minimal Cuntz-Krieger Dilations and
 Representations of Cuntz-Krieger Algebras}}

\end{center}

\vsp \vsp \vsp

\begin{center}
{\sc B. V. Rajarama Bhat,  Santanu Dey and  Joachim Zacharias }
\vsp \vsp

{\bf 24 May 2005}
\end{center}

\vsp \vsp \vsp \vsp

\begin{center}
{\underline {\bf Abstract}}
\end{center}

\vsp
Given a contractive tuple of Hilbert space operators satisfying certain $A$-relations
we show that there exists a unique minimal dilation to generators of Cuntz-Krieger
algebras or its extension by compact operators. This Cuntz-Krieger dilation
can be obtained from the classical minimal isometric dilation as a certain
maximal $A$-relation piece. We define a maximal piece more generally for
a finite set of polynomials in $n$ noncommuting variables.
We classify all representations of Cuntz-Krieger algebras $\CO_A$
obtained from dilations of commuting tuples satisfying
$A$-relations. The universal properties of the minimal
Cuntz-Krieger dilation and the WOT-closed algebra generated
by it is studied in terms of invariant subspaces.

\vsp \vfill
----------------------------------------------------------------------



\newpage

\begin{section}{Introduction}
\setcounter{equation}{0}
Cuntz-Krieger algebras were introduced by
J. Cuntz and W. Krieger in [CK] as examples of simple purely infinite $C^*$-algebras
not stably isomorphic to Cuntz algebras. Let
$A=(a_{ij})_{n\times n}$ be a square $0-1$-matrix i.e. $a_{ij} \in \{ 0 , 1\}$ and each
row and column has at least one non-zero entry. The Cuntz-Krieger algebra
$\CO_A$ is defined as follows:

\begin{Definition} {\em $\CO_A$ is the universal
$C^*$-algebra generated by $n$ partial isometries $s_1,
\cdots,s_n$ with orthogonal ranges satisfying
\begin{equation}
\begin{array}{c}
s^*_is_i= \sum^n_{j=1} a_{ij} s_js^*_j\\
\\
I=\sum^n_{i=1}s_is^*_i.
\end{array}
\end{equation}
We denote the tuple $(s_1,\cdots, s_n)$ by $\us.$}
\end{Definition}

Notice that $s_is_j = a_{ij}s_is_j$ for all $s_i$ and $s_j$ in $\CO_A$.
In this paper we study dilations related to these algebras.  Equations (1.1)
are called {\em Cuntz-Krieger relations.} An $n$-tuple of
bounded operators $\uT=(T_1,\cdots,T_n)$ on a Hilbert space is
said to be a {\em contractive $n$-tuple} if $T_1T^*_1+ \cdots + T_nT^*_n
\leq I.$ Such tuples are also called row contractions as the
condition is equivalent to saying that the operator $(T_1,\cdots,T_n)$
from $\CH \oplus \cdots \oplus \CH$ ($n$-times) to $\CH$ is a
contraction. We will only consider contractive tuples.
For such tuples Davis [D], Bunce [Bu],
Frazho [Fr] and more extensively Popescu ([Po1-6], [AP]) constructed
dilations consisting of isometries with orthogonal ranges. Under natural
minimality conditions this dilation is unique up to unitary equivalence.
We refer to it as the {\em minimal isometric dilation} or the
{\em standard (noncommuting) dilation} $\hat{\uT}$ acting on
$\hat{\CH}$ (c.f.\;[BBD]).
In case $T_1T^*_1+ \cdots + T_nT^*_n= I$ we will sometimes call $\uT$ {\em unital}.
$\uT$ is unital iff its standard dilation is unital.
The standard noncommuting dilation has similar characterizations as classical dilations of
single operators notably the minimal normal extension of subnormal operators.
The universal role of the unilateral shift for single operators is played
by the tuple of creation operators on the full Fock space.

However, tuples are more complex than single operators and one may impose
symmetry conditions on the tuple and study dilations within this restricted class
of tuples.
For instance if all operators in the tuple commute i.e.\;$\uT$ is a
commuting tuple,  Arveson [Ar3] showed that there is a
unique minimal commuting dilation with similar properties.
Crucial in his approach is the tuple of creation operators on symmetric Fock space
playing the role of the shift for single operators.
In [BBD] the relation between the minimal commuting and the standard dilation has
been investigated. It was shown that
for every contractive tuple there is a maximal
subspace on which it forms a commuting tuple
and that the minimal commuting dilation is precisely
this maximal commuting piece of the standard noncommuting
dilation.

In this article we consider the
following class of tuples and investigate dilations within this class
and their connection with the commuting and noncommuting standard dilation.

\begin{Definition}
{\em Let $A = (a_{ij})_{n\times n}$ be a 0-1-matrix. A contractive $n$-tuple
$\uT$ is an {\em $A$-relation tuple\/} or is said to satisfy {\em $A$-relations \/} if
$T_iT_j=a_{ij}T_iT_j$ for $1 \leq i,j \leq n$.}
\end{Definition}

Given such a tuple $\uT$ there is a unique minimal dilation to partial
isometries satisfying $A$-relations i.e.\;generators of a Cuntz-Krieger algebras if $\uT$
is unital or an extension of it by compact operators if $\uT$ is contractive.
We call this the minimal Cuntz-Krieger dilation of $\uT$.
It will be denoted by $\tilde{\uT}$ and acts on $\tilde{\CH}$.

For an arbitrary tuple we define a maximal $A$-relation
piece and compare the maximal $A$-relation piece of
the standard isometric dilation with the minimal Cuntz-Krieger
dilation. As for commuting tuples both turn out to be the same.
We also prove similar results for the maximal commuting
$A$-relation piece.

We begin in section 2 by defining the maximal piece of a tuple of operators
with respect to a finite set of polynomials in $n$-noncommuting variables.
Similarly to results by Arias and Popescu [AP] there is a canonical homomorphism
between the WOT-closed (non-selfadjoint) algebra generated by them and the WOT-closed
algebra generated by the original tuple modulo a two-sided ideal.
The maximal $A$-relation piece of a tuple of $n$-isometries with orthogonal
ranges is a special case of this. In fact the maximal commuting piece and
maximal $q$-commuting piece (c.f.\:[BBD], [De], [AP]) can all be treated using
this approach.

In section 3 we show that the maximal $A$-relation piece of
the standard dilation of an $n$-tuple satisfying $A$-relations is
the minimal Cuntz-Krieger dilation. The section begins with two different constructions
of minimal Cuntz-Krieger dilations, one using positive definite kernels the other
using a modification of Popescu's Poisson transform.

In section 4 we study the minimal Cuntz-Krieger dilation of commuting $A$-relation
tuples from a representation theoretical point of view. If such a tuple is also unital
then it determines a unique representation of the Cuntz-Krieger algebra $\CO_A$.
Generalizing results from [BBD] for $\CO_n$ we are able to show that these
representations are determined by the GNS-representations of analogues of Cuntz states.

Based on ideas of Bunce, Popescu [Po5] showed that the minimal isometric dilation
can be characterized by a universal property of the $C^*$-algebra generated by it.
In section 5 we first point out that minimal Cuntz-Krieger dilations can
be characterized in a similar way. We study the structure of the WOT-closed
algebra generated by the operators constituting the minimal Cuntz-Krieger dilation and
describe this algebra by making use of its invariant and wandering
subspaces. We use techniques of Davidson, Kribs and Pitts ([DPS], [DP2]) to understand
the structure of `free semigroup algebras' i.e.\;WOT-closed algebras generated by a
finite number of isometries with orthogonal ranges. In this section many proofs are only
sketched or omitted.

All Hilbert spaces in this paper are complex and separable. We denote
the full Fock space over $\CL$ by $\Gamma(\CL)$ which is defined as
$$
\Gamma (\CL) = \C \oplus \CL \oplus  \CL^{\otimes ^2} \oplus
\CL ^{\otimes ^3} \oplus \ldots.
$$
Let the vacuum vector $1 \oplus 0 \oplus 0 \oplus  \ldots $
be denoted by $\omega.$
$\C^n$ is the $n$-dimensional complex Euclidean space with
standard orthonormal basis $\{e_1,\cdots, e_n\}$. The left
creation operator $L_i$ on
$\Gamma(\C^n)$
is defined by
$$
L_ix= e_i \otimes x,
$$
where $1 \leq i \leq n$ and $x \in \Gamma(\C^n).$  The $L_i$'s are clearly
isometries with orthogonal ranges. We denote the
tuple $(L_1,\cdots,L_n)$  by $\uL$. Also $\sum_i L_iL^*_i=I-P_0$
where $P_0$ is the projection onto the vacuum space.

Let $\Lambda$ be the set $\{ 1,2,\cdots, n\}$ and $\Lambda^m$
the $m$-fold cartesian product of $\Lambda$ for $m \in \N$.
For an operator tuple $(T_1,\cdots,T_n)$ on a Hilbert space $\CH$ and for
$\alpha=(\alpha_1, \alpha_2 ,\ldots,\alpha_m)= \alpha_1 \alpha_2 \cdots \alpha_m$
in $\Lambda^m$, the operator $T_{\alpha_1}T_{\alpha_2}\cdots T_{\alpha_m}$ will be
denoted by $\uT^{\alpha}$. Let $\tilde{\Lambda}$ denote
$\cup^\infty_{m=0} \Lambda^m$, where $\Lambda^0$ is $\{0\}$ and
$\uT^0$ is the identity operator. We may think of elements in $\tilde{\Lambda}$
as words with concatenation as  product written $\alpha \beta$.
$\tilde{\Lambda}$ is the free semigroup with $n$ generators. Given a 0-1-matrix $A$
as above we can define
$\Lambda_A^m=\{ \alpha_1 \alpha_2 \cdots \alpha_m : a_{ \alpha_i , \alpha_{i+1}} = 1
\text{ for } i=1 , \ldots ,n-1 \}$ and
the subsemigroup $\tilde{\Lambda}_A=\bigcup^{\infty}_{m=0} \Lambda_A^m$.
$o(\alpha)$ and $t(\alpha)$ denote the first and last letter (i.e.\;index) of $\alpha$.

\begin{Definition}
{\em Let $\CH$ and $\CL$ be two Hilbert spaces such that $\CH$ is
a closed subspace of $\CL$ and $\uT,~~\uR$ be $n$-tuples of
operators on $\CH,~~\CL$ respectively. Then $\uR$ is a {\em
dilation \/} of $\uT$ or $\uT$ a {\em piece}
of $\uR$ if
$$
R^*_ih=T^*_ih
$$
for all $h  \in \CH, ~~1\leq i \leq n$.   A dilation is said to be {\em minimal dilation \/} if
$$
\overline{\mbox{span}}\{ \uR^{\alpha}h: \alpha \in \tilde{\Lambda}, h \in \CH\}=\CL.
$$
\begin{enumerate}
\item A dilation $\uR$ of $\uT$ is said to be {\em
isometric \/} if $\uR$ consists of isometries with orthogonal ranges.
\item When $\uT$ satisfies $A$-relations, a dilation $\uR$ of $\uT$
is said to be a {\em Cuntz-Krieger dilation \/} if $\uR$ consists
of partial isometries with orthogonal ranges satisfying $A$-relations and
\begin{equation}
R^*_iR_i=I-\sum^n_{j=1}(1-a_{ij})R_jR^*_j=P_0 + \sum_{j=1}^n a_{ij} R_jR_j^*,
\end{equation}
where $P_0 = 1-\sum_{j=1}^n  R_jR_j^*$.
\end{enumerate} }
\end{Definition}

Thus, like for isometric dilations,
$C^*(\uR)= \overline{\text{span}}\{ \uR^{\alpha} (\uR^{\beta})^* : \alpha , \beta \in \tilde{\Lambda} \}$
for any Cuntz-Krieger dilation $\uR$ of $\uT$. Moreover, since
\begin{equation}
(\uR^{\alpha})^* \uR^{\beta} = \delta_{\alpha,\beta} R^*_{t(\alpha)} R_{t(\alpha)}= \delta_{\alpha,\beta} \Big(I-\sum^n_{j=1}(1-a_{t(\alpha),j})R_jR^*_j \Big)
\end{equation}
whenever  $\alpha \neq 0$, it follows that if $\uR_1$ and $\uR_2$
are any two minimal Cuntz-Krieger dilations then
$\sum \uR_1^{\alpha_i} h_i \mapsto \sum \uR_2^{\alpha_i} h_i$ extends to a unitary equivalence.

Given any dilation $\uR$ of $\uT$ all $R_i^*$ leave $\CH$ invariant and
if $p,q$ are polynomials in $n$-noncommuting variables then
$$
\uT^\alpha(\uT^\beta)^*=P_{\CH}\uR^\alpha (\uR ^\beta)^*|_{\CH}
\mbox{~~~~~and~~~~~}
p(\uT)(q(\uT))^*=P_{\CH}p(\uR)(q(\uR))^*|_{\CH.}
$$
It follows that if $\uR$ is a Cuntz-Krieger dilation, then there is a unique completely positive map
$\rho :C^*(\uR) \to C^*(\uT)$ mapping $\uR^{\alpha} (\uR^{\beta})^*$ to $\uT^{\alpha} (\uT^{\beta})^*$.

\smallskip

Finally let us recall a concept needed later. For an $n$-tuple $\uR$ of bounded operators on $\CL$, a subspace $\CK$ of $\CL$ is said to be {\em wandering \/} for the tuple
if $\uR^\alpha
\CK$ are pairwise orthogonal for all $\alpha \in \tilde{\Lambda}$.
\end{section}

\begin{section}{Maximal $A$-relation Piece and $A$-Fock Space}
\setcounter{equation}{0}

We begin with an $n$-tuple of bounded operators $\uR$ on a Hilbert
space $\CL$ and a finite set of polynomials $\{p_\xi\}_{\xi \in \CI}$ in $n$-noncommuting
variables with finite index set $\CI$. Consider
$$
{\mathcal C}(\uR )=\{ \CM :  R_i^* \CM \subseteq \CM \text{ and } (p_\xi(\uR))^*h=0, \forall h\in \CM,
1\leq i\leq n, \xi \in \CI \}.
$$
${\mathcal C}(\uR)$ consists of all co-invariant subspaces of $\uR $ such
that the compressions form a tuple $\uR ^p=(R_1^p, \ldots , R_n^p)$ satisfying $p_\xi
(\uR^p)=0$ for all $\xi \in \CI.$ It is a complete lattice, in the sense
that arbitrary intersections and closed spans of arbitrary unions of
such spaces are again in this collection. Its maximal
element is denoted by $\CL^p (\uR)$ (or by $\CL^p $ when the
tuple under consideration is clear).
Since $(p_\xi (\uR))^*(\uR ^{\alpha})^* \CM = 0$ for all $\CM \in
{\mathcal C}(\uR )$,
$\alpha \in \tilde{\Lambda}$ and $\xi \in \CI$ we have
$\CL^p (\uR) \subseteq \bigcap_{\alpha \in \tilde{\Lambda} , \xi \in \CI}
\ker (p_\xi (\uR)^* (\uR^{\alpha })^*)$. On the other hand this intersection lies
in ${\mathcal C}(\uR )$ hence
$$
\CL^p (\uR) = \bigcap_{\alpha \in \tilde{\Lambda} , \xi \in \CI} \ker (p_\xi (\uR)^* (\uR^{\alpha})^*) =
\left[ \bigvee_{\alpha \in \tilde{\Lambda}, \xi \in \CI} \uR^{\alpha} p_\xi (\uR)(\CL) \right]^{\perp}.
$$
Therefore we have:
\begin{Lemma}
Let $\uR$ be an $n$-tuple of operators on a Hilbert space
${\mathcal L}$ and
$\CK = \overline{span}\{ \uR ^{\alpha } p_\xi (\uR) \uR^\beta h: h\in \CL, \; \xi \in \CI \mbox{~and~} \alpha, \beta \in \tilde
{\Lambda }\}$. Then $\CL ^p(\uR) = \CK ^{\perp } =\{ h\in \CL: (\uR ^\alpha p_\xi(\uR)\uR^\beta)^* h=0, \forall \xi
\in \CI \mbox{~and~} \alpha, \beta \in {\tilde \Lambda}\}.$
\end{Lemma}

\noindent
$\CL ^p(\uR)$ can also be thought of as follows.
Let $\CR$ be the (non self-adjoint) WOT-closed algebra generated by $\uR$. Then
$\CL ^p(\uR)=( \CJ \CL)^{\perp}$, where
$$
\CJ=\overline{\mbox{span}}^w \{\uR^\alpha p_{\xi}(\uR)\uR^\beta:
\alpha, \beta \in \tilde{\Lambda}, \xi \in \CI \} \subseteq \CR
$$
is the WOT-closed ideal generated by $\{ p_{\xi}(\uR) : \xi \in \CI \}$.

\begin{Definition}{\em   The {\em maximal  piece\/ } of
$\uR$ with respect to $\{p_{\xi}\}_{\xi \in \CI}$ is defined as
the piece obtained by compressing
$\uR$ to the maximal element $\CL^p (\uR)$ of ${\mathcal C}(\uR)$ denoted
by $\uR ^p=(R_1^p, \ldots , R_n^p)$. The maximal piece is said to be {\em
trivial \/} if the space $\CL^p(\uR)$ is the zero space.}
\end{Definition}

Let $\CR^p$ be the  WOT-closed algebra generated by $\uR^p$.
By co-invariance the map
$$
\Psi_{\CL ^p(\uR)}: \CR \to \CR^p ,\;\;\;\; X \mapsto P_{\CL ^p(\uR)} X=P_{\CL ^p(\uR)} XP_{\CL ^p(\uR)}
$$
is a WOT-continuous homomorphism of $\CR$ whose kernel is a WOT-closed ideal of $\CR$.
Since $\Psi_{\CL^p} ( p_{\xi}(\uR)) = p_{\xi}(\uR^p)=0$ for all $\xi \in \CI$ we
certainly have $p_{\xi}(\uR) \in \ker \Psi_{\CL^p}$ i.e. $\CJ \subseteq \ker \Psi_{\CL^p}$.
Thus there is a canonical surjective and contractive homomorphism
$$
\Psi : \CR /\CJ \to \CR^p.
$$

When the polynomials are $p_{(l,m)}=z_lz_m-a_{lm}z_lz_m, (l,m) \in
\{1,\cdots, n\}\times \{1, \cdots, n\}=\CI$ we call $\CL^p(\uR)$
the {\em maximal $A$-relation subspace} and the corresponding piece the
{\em maximal A-relation piece} $\uR^A$. The  maximal $A$-relation subspace is explicitly given by
$$
\CL_A(\uR)=\{  \uR^{\alpha}(R_i R_j - a_{ij}R_i R_j) h : h \in \CL , \; \alpha \in \Lambda, \; i,j=1,\ldots,n\}^{\perp}.
$$
When the noncommuting polynomials are $q_{(l,m)} = z_lz_m - z_m z_l$ then we obtain
the maximal commuting subspace
$$
\CL^c(\uR)=\{  \uR^{\alpha}(R_i R_j - R_j R_i) h : h \in \CL , \; \alpha \in \Lambda, \; i,j=1,\ldots,n\}^{\perp}.
$$
studied in [BBD].
\begin{Lemma}
Let $\uR$ and $\uT$ be two $n$-tuples of bounded operators on $\CM$ and $\CH$ respectively.
\begin{enumerate}
\item The maximal $A$-relation piece of $(R_1\oplus T_1,\cdots,R_n \oplus T_n)$
is $(R^A_1 \oplus T^A_1,\cdots,R^A_n\oplus T^A_n)$ acting on
$\CM_A \oplus \CH_A$ and the maximal $A$-relation piece of
$(R_1\otimes I,\cdots, R_n\otimes I)$ acting on $\CM \otimes \CH$
is $(R^A_1\otimes I,\cdots, R^A_n\otimes I)$ on $\CM_A
\otimes\CH.$
\item Suppose $\H \subseteq \CM$ and $\uR$ is a dilation of $\uT$
then $\uR^A$ is a dilation of $\uT^A$ with
$\CH_A(\uT)=\CM_A(\uR)\cap \CH.$
\end{enumerate}
\end{Lemma}
\noindent{Proof:}
Follows from Lemma 4 (compare with [BBD] for part (2)).
\qed

\medskip

Cuntz-Krieger relations are naturally related to $A$-Fock space, a variant of the usual
Fock space.

\begin{Definition}
{\em For a given $A=(a_{ij})_{n\times n}$ as above, the {\em
$A$-Fock space \/} is defined as the maximal $A$-relation subspace
$(\Gamma(\C^n))_A(\uL)$ with respect to the left creation
operators. It is denoted by $\Gamma_A$. We also define
the $n$-tuple $\uS=(S_1,\cdots,S_n)$, where the $S_i$'s are the
compressions of left creation operators $L_i$ onto
$\Gamma_A.$}
\end{Definition}
\noindent
The $A$-Fock space has a very concrete description justifying the terminology.

\begin{Proposition}
$\Gamma_A=\overline{\mbox{span}}\{ e^{\alpha}:
\alpha \in \tilde{\Lambda}_A \}$ and $\uL^A = \uS$.
\end{Proposition}
\noindent {\sc Proof:} Let $\alpha \in {\Lambda}^m$ be such that
there exist $1\leq k \leq m-1$ for which $a_{\alpha_k
\alpha_{k+1}}=0.$ Denoting $\alpha_{k}, \alpha_{k+1}$ by $s, t,$
it is clear that
$$
e^\alpha  \in \overline{\mbox{span}}\{ \uL^\gamma
(L_sL_t-a_{st}L_sL_t)h: h \in \Gamma(\C^n), \gamma \in \tilde{\Lambda}
\},
$$
which implies that such $e^{\alpha}$ are orthogonal to
$\Gamma (\C^n)_A(\uL)$, whereas if for all $1 \leq k \leq m-1, ~~
a_{\alpha_k\alpha_{k+1}}=1$ then for all $1\leq i,j \leq m-1,
\beta \in \tilde{\Lambda}, h \in \Gamma(\C^n)$
$$
\langle  e^\alpha, \uL^\beta (L_iL_j-a_{ij}L_iL_j)h \rangle =0,
$$
and thus such $e^{\alpha} \in \Gamma_A$. By taking
completions the Proposition follows.
\qed

\medskip

Similar Fock spaces were also studied by Muhly [Mu], Solel and others.

Now suppose that $\alpha = \alpha_1 \cdots \alpha_m \in \Lambda_A^m$ and $m > 0$, then
$$
S_i e^\alpha = P_{\Gamma_A} L_ie^\alpha= \left\{
\begin{array}{cc}
e_i & \mbox{~if~} |\alpha|=0 \\
a_{i \alpha_1} e_i \otimes e^\alpha & \mbox{~if~} |\alpha| \geq 1
\end{array} \right.$$
$$
S^*_i e^{\alpha}=  L^*_i e^\alpha=\left\{
\begin{array}{cc}
0 & \mbox{~if~} |\alpha|=0 \\
\delta_{i \alpha_1} \omega & \mbox{~if~} |\alpha|=1 \\
\delta_{i\alpha_1}e_{\alpha_2} \otimes \cdots \otimes e_{\alpha_m} & \mbox{~if~}
|\alpha| > 1,
\end{array} \right.
$$
$$ S^*_i S_ie^\alpha = \left\{ \begin{array}{cc} \omega & \mbox{~if~}
|\alpha|=0 \\
a_{i \alpha_1} e^\alpha & \mbox{~if~} |\alpha| \geq 1
\end{array} \right.
\mbox{~and~}
S_iS^*_i e^{\alpha}= \left\{ \begin{array}{cc} 0 & \mbox{~if~} |\alpha|=0 \\
\delta_{i \alpha_1} e^\alpha & \mbox{~if~} |\alpha|=1.
\end{array} \right.
$$

\begin{Proposition}
The maximal $A$-relation piece of an $n$-tuple of isometries with
orthogonal ranges is an $n$-tuple of partial isometries with orthogonal
ranges.
\end{Proposition}
\noindent {\sc Proof:} Let
$\uV=(V_1,\cdots,V_n)$ be an $n$-tuple of
isometries with orthogonal ranges on a Hilbert space $\CK.$ Fix a
matrix $A=(a_{ij})_{n\times n}$ as above and denote the projection
onto $\CK_A(\uV)$ by $P.$ Any $k_A$ in $\CK_A(\uV)$ can be
written as $k_A= \oplus^n_{p=1} V_p k_p \oplus k_0$ for some
$k_p\in \CK, ~~1\leq p \leq n$ and some $k_0 \in
(I-\sum^n_{p=1}V_pV^*_p)\CK$. (Any $k \in \CK$ can be written in this form.)
Clearly $k_0$ is in $\CK_A(\uV)$ using Lemma 4, since the ranges of the $V_q$'s and
$I-\sum_i V_i V^*_i$ are all mutually orthogonal. Similarly one
observes that the other $k_p$'s also
belong to $\CK_A(\uV)$ as for $k \in \CK, \alpha \in
\tilde{\Lambda}$
\begin{eqnarray*}
\langle
k_p,\uV^\alpha(V_i V_j -a_{ij} V_i V_j)k\rangle
&=& \langle V_p
k_p ,V_p \uV^\alpha(V_i V_j -a_{ij}
V_i V_j)k \rangle\\
&=& \langle \oplus^n_{q=1} V_q k_q \oplus k_0,
V_p \uV^\alpha(V_i V_j-a_{ij}
V_i V_j)k\rangle\\
&=&\langle
k_A,V_p \uV^\alpha(V_i V_j-a_{ij}
V_i V_j)k\rangle=0,
\end{eqnarray*}
where again we use that the ranges of the $V_q$'s and
$I-\sum_i V_i V^*_i$ are all mutually orthogonal. Next
we show that
\begin{equation}
PV_i k_0= V_i k_0,
\end{equation}
and
\begin{equation}
PV_i V_p k_p= a_{ip} V_i V_p k_p.
\end{equation}
Equation (2.1) follows from $\langle V_i k_0, \uV^\beta
(V_s V_t -a_{st} V_s V_t) k \rangle=0,$ for all
$\beta \in \tilde{\Lambda},1 \leq s,t \leq n, ~~k \in \CK$
(since $k_0$ is orthogonal to the range of $V_t, 1 \leq t \leq n$).
 When $a_{ip}=0,$ we have $ PV_i V_p k_p=P(V_i V_p -a_{ip}
 V_i V_p ) k_p=0=a_{ip}
 V_i V_p k_p.$ So it is enough to show for
$a_{ip}=1$ that
$V_i V_p k_p \in \CK_A(\uV).$ When $|\alpha|
> 1$ or $|\alpha|=0,$ it is easy to see that for $1 \leq s,t \leq n, ~~k \in \CK$
 \begin{equation}
  \langle V_iV_pk_p, \uV^\alpha (V_s
 V_t-a_{st}V_sV_t)k \rangle =0
 \end{equation}
 as the $V_i$'s are isometries with orthogonal ranges and $k_p \in \CK_A(\uV).$
 When $|\alpha|=1,$
\begin{eqnarray*}
\langle V_iV_pk_p, V_i (V_p
 V_t-a_{pt}V_pV_t)k \rangle
 &=& \langle V_pk_p,  (V_p
 V_t-a_{pt}V_pV_t)k \rangle\\
 &=&\langle \oplus^n_{s=1}V_s k_s\oplus k_0,  (V_p
 V_t-a_{pt}V_pV_t)k \rangle=0.
\end{eqnarray*}
Clearly equation (2.3) holds in all other cases when
$| \alpha |=1.$ So equation (2.2) is true in general and we have
\begin{eqnarray*}
V^A_i (V^A)^*_i V_i ^A k_A &= & P V_iV^*_i  P
V_i k_A\\
&=& P V_iV^*_i P(\oplus_p V_iV_p k_p \oplus
V_ik_0)\\
& =&P V_iV^*_i(\oplus^n_{p=1} a_{ip}
V_iV_pk_p \oplus V_i k_0)\\
&=& \oplus^n_{p=1} a_{ip} P V_iV_pk_p \oplus P
V_i k_0 \\
&=& PV_i k_A= V^A_i k_A.
\end{eqnarray*}
Thus the $V^A_i$'s are partial isometries. Now the assertion of the
Proposition
that for $1 \leq i \neq j \leq n$, the range of $V^A_i$ is orthogonal
to the range of $V^A_j$ can be proved in the following way:
\begin{eqnarray*}
(V^A_j)^*V^A_ik_A& = & V^*_j  P V_i k_A\\
&=&  V^*_j P V_i (\oplus^n_{p=1} V_p k_p \oplus
k_0)\\
&=& V^*_j (\oplus_p a_{ip}V_iV_p k_p \oplus
V_ik_0)=0.\\
\end{eqnarray*}
Alternatively this follows since $\uV^A$ is contractive.
\qed

\begin{Corollary}
The following holds for $\uS:$
\begin{enumerate}
\item $I-\sum^n_{i=1}S_iS^*_i=P_0$ where $P_0$ is the projection
onto the vacuum space.
\item The $S_i$'s are partial isometries with orthogonal ranges.
\item $S^*_iS_i=I-\sum^n_{j=1}(1-a_{ij})S_jS^*_j= P_0 + \sum^n_{j=1}a_{ij}S_jS^*_j$
\end{enumerate}
\end{Corollary}
\noindent{\sc Proof:}
\begin{enumerate}
\item $I-\sum S_iS^*_i=P_{\Gamma_A}(I-\sum L_iL^*_i)
P_{\Gamma_A}=P_0$.
\item Follows from Proposition 9 or can be checked directly from the relations before Proposition 9.
\item Suppose $e^\alpha \in \Gamma_A$, and when $ |\alpha| > 0$ let
$\alpha=\alpha_1 \cdots \alpha_m$. Then
$$[I-\sum_j(1-a_{ij})S_jS^*_j]e^\alpha= \left\{ \begin{array}{cc} \omega
& \mbox{~if~}
|\alpha|=0 \\
a_{i \alpha_1} e^\alpha & \mbox{~if~} |\alpha| \geq 1.
\end{array} \right.$$
\end{enumerate} \qed
\end{section}

\begin{section}{Minimal Cuntz-Krieger Dilations and Standard
Noncommuting Dilations}
\setcounter{equation}{0}

\noindent
The main purpose of this section is to prove the following result.

\begin{Theorem}
Let $\uT$ be a contractive $A$-relation tuple on $\CH$. Then there exists a minimal
Cuntz-Krieger dilation $\tilde{\uT}$ on $\tilde{\CH}$ unique up to unitary equivalence.
$\tilde{\uT}$ is unital iff $\uT$ is unital.
\end{Theorem}
Uniqueness has already been pointed out and the last part follows since
$\sum_{|\alpha|=k} \tilde{\uT}^\alpha (\tilde{\uT}^\alpha )^*$ form a decreasing sequence
of projections converging weakly to a limit projection $P$. If $\uT$ is unital then $P_{\CH} \leq P$
and so $P=1$ by minimality.

We will give two proofs of the existence. The first is direct and uses positive
definite kernels. It gives an explicit construction of the dilation Hilbert space and the dilated tuple.
The second construction is an adaptation of Popescu's Poisson transform method which uses
completely positive maps. Though elegant it is less direct. Then we show that the minimal
Cuntz-Krieger dilation can also be obtained as the maximal $A$-relation piece of the standard
dilation.

\medskip

\noindent
\textsc{First proof via positive definite Kernels}

\medskip

Let $\uT=(T_1,\cdots,T_n)$
be a contractive $n$-tuple on a Hilbert space $\CH$ satisfying $A$-relations.
Assume we already found the minimal Cuntz-Krieger dilation $\tilde{\uT}$ on the Hilbert space $\tilde{\CH}$.
Then
$$
K((\alpha,u),(\beta,v)):=\langle \tilde{\uT}^{\alpha} u ,  \tilde{\uT}^{\beta} v \rangle
$$
clearly defines a positive definite kernel on the set
$$
X= \tilde{\Lambda}_A \times \CH.
$$
By minimality $\tilde{\CH} = \overline{\text{span}}\{\tilde{\uT}^{\alpha} u : u \in \CH, \alpha \in \tilde{\Lambda}_A \}$
which is precisely the kernel Hilbert space. Moreover $\tilde{T}_i$ corresponds to the map $(\alpha , u) \mapsto (i \alpha , u)$.

Using co-invariance of $\CH$ under $\tilde{\uT}$ and the relation $\tilde{T}_i^* \tilde{T}_i = I - \sum_j (1-a_{ij}) \tilde{T}_j \tilde{T}_j^* $
we find that
$$
K((\alpha,u),(\beta,v))
=
\left\{
\begin{array}{cc}
\langle u, v\rangle &\mbox{if~}\alpha=\beta=0\\
\langle u, [I-\sum_k
(1-a_{t(\alpha)k}) T_kT^*_k]
v\rangle & \mbox{if~} \alpha=\beta\neq 0 \\
\langle u , \uT^\gamma v
\rangle
&\mbox{if~} \beta=\alpha \gamma  \\
\langle u, (\uT^\gamma)^* v
\rangle & \mbox{if~}\alpha = \beta \gamma \\
0 & \mbox{otherwise.}
\end{array}
\right.
$$
and this kernel depends only on $\uT$. We will show directly by induction that
the kernel $K$ thus defined is always positive definite.

To simplify the calculations we assume that $\uT$ is unital i.e. $\sum_i T_i T_i^* =I$.
There is no loss in doing so since positivity for the kernel defined by the $(n+1)$-tuple
$(T_1, \ldots ,T_n,(I-\sum_i T_iT_i^*)^{1/2})$ implies positivity of $K$.
Under this assumption $T_iQ_i= T_i$, where $Q_i=I-\sum_k(1-a_{i k}) T_kT^*_k = \sum_{j=1}^n a_{i,j} T_j T_j^*$
and $K((\alpha,u),(\alpha,v))= \langle u , Q_{t(\alpha)} v \rangle$, whenever $\alpha \neq 0$.

Now let $A^{(m)}$ denote operator matrices with entries in $B (\CH)$ indexed by $\alpha,\beta \in \tilde{\Lambda},$ where $|\alpha|,|\beta| \leq m$ and define $K^{(m)}=(K^{(m)}_{\alpha,\beta})$ by
$$
K^{(m)}_{\alpha,\beta}:=\left\{
\begin{array}{cc} I &\mbox{if~}\alpha=\beta=0 \\
Q_{t (\alpha )}
&\mbox{if~} \alpha=\beta\neq 0 \\
\uT^\gamma &\mbox{if~} \beta=\alpha \gamma \\
(\uT^\gamma)^* & \mbox{if~} \alpha = \beta \gamma \\
0& \mbox{otherwise.}
\end{array} \right.
$$
i.e. $K^{(m)}$ is a compression of $K$. Clearly it suffices to show that all $K^{(m)}$ are positive.
For $m=1$ this follows form the equation
$$
\left[
\begin{array}{ccccc}
I      &  T_1  &  T_2  & \cdots & T_n \\
T_1^*  &  Q_1  &   0   & \cdots &  0  \\
T_2^*  &   0   &  Q_2  & \cdots &  0  \\
\vdots &       &       &        & \vdots \\
T_n^*  &   0   &       & \cdots       & Q_n
\end{array}
\right]
=
$$
$$
=
\left[
\begin{array}{ccccc}
I      &  T_1  &  T_2  & \cdots & T_n \\
0      &   I   &   0   & \cdots &  0  \\
0      &   0   &   I   & \cdots &  0  \\
\vdots &       &       &        & \vdots \\
0      &   0   &       &  \cdots      &  I
\end{array}
\right]
\left[
\begin{array}{ccccc}
 0     &  0    &   0   & \cdots &  0  \\
 0     &  Q_1  &   0   & \cdots &  0  \\
 0     &   0   &  Q_2  & \cdots &  0  \\
 0     &       &       &        & \vdots \\
 0     &   0   &       & \cdots       & Q_n
\end{array}
\right]
\left[
\begin{array}{ccccc}
I      &   0   &   0   & \cdots &  0  \\
T_1^*  &   I   &   0   & \cdots &  0  \\
T_2^*  &   0   &   I   & \cdots &  0  \\
\vdots &       &       &        & \vdots \\
T_n^*  &   0   &       & \cdots       &  I
\end{array}
\right]
$$
If $m>1$ define matrices $L_1, L_2,\ldots ,L_{m-1}$ by
\[
L_{k; \alpha,\beta}:=
\left\{
\begin{array}{cc} T_i & \mbox{if~} \beta = \alpha i ,|\beta|=k \\
I & \mbox{if~} \alpha=\beta \text{ and } |\beta | \geq k \\
0 & \mbox{otherwise.}
\end{array}
\right.
\]
Finally let
\[
Q_{\alpha,\beta}^{(m)}:=
\left\{
\begin{array}{cc} Q_{t({\alpha})} &
\mbox{if ~}\alpha=\beta \text{ and } |\alpha |=m \\
0 & \mbox{otherwise.}
\end{array}
\right.
\]
Then it is not hard to check that
$$
K^{(m)} = L_1 L_2 \ldots L_{m-1} Q^{(m)} L_{m-1}^* \ldots L_2^* L_1^*
$$
which shows that $K$ is positive. By Kolmogorov's Theorem there exists
a Hilbert space $\tilde{\CH}$ and an injective map $\lambda: X \to \tilde{\CH}$
such that $\overline{\mbox{span}}\{\lambda (\alpha,u): 1 \leq i
\leq n, \alpha \in \tilde{\Lambda}, u \in \CH\}=\tilde{\CH}$ and
$$
K((\alpha,u) ,(\beta,v))=\langle \lambda(\alpha,u) ,\lambda(\beta,v) \rangle.
$$
It remains to show that
$\tilde{\uT}=(\tilde{T}_1,\cdots,\tilde{T}_n)$ consisting of maps $\tilde{T}_i:\tilde{\CH} \to \tilde{\CH}$ defined by
$$
\tilde{T}_i\lambda(\alpha,u)=
\lambda( i \alpha,u),
$$
constitute a tuple $\tilde{\uT}$ which is the minimal Cuntz-Krieger dilation of $\uT$.
First note that $\tilde{T}_i$ is a well-defined contraction.
Indeed, thinking of $K$ as a block matrix we have
$$
K_{i\alpha,i\beta}=\left\{
\begin{array}{cc}  Q_i &\mbox{if~}\alpha=\beta=0 \\
a_{io(\alpha)}a_{io(\beta)} K_{\alpha , \beta}
& \mbox{otherwise,}
\end{array} \right.
$$
where we define $a_{io(\alpha)}=1$ for all $i$ if $\alpha=0$.
So if $F \subseteq X$ is a finite subset then
\begin{eqnarray*}
\Big\| \sum_{(\alpha,u) \in F} \lambda (i\alpha,u) \Big\|^2
&=&
\Big\| \sum_{(\alpha,u) \in F} a_{io(\alpha)} \lambda (i\alpha,u) \Big\|^2 \\
&=&
\sum_{\alpha = \beta=0} \langle u ,Q_i v \rangle +
\sum_{\alpha \neq 0 \text{ or }  \beta \neq 0} a_{io(\alpha)} a_{io(\beta)}\langle u ,K_{\alpha , \beta} v \rangle \\
& \leq &
\sum_{\alpha = \beta=0} \langle u ,v \rangle +
\sum_{\alpha \neq 0 \text{ or }  \beta \neq 0}  \langle u ,K_{\alpha , \beta} v \rangle \\
&=&
\Big\| \sum_{(\alpha,u) \in F} \lambda (\alpha,u) \Big\|^2.
\end{eqnarray*}
For $i \neq j$
\begin{eqnarray*}
\langle \tilde{T}_i\lambda(\alpha ,u),
\tilde{T}_j\lambda(\beta,v)\rangle
&=& \langle \lambda(i\alpha,u),\lambda(j\beta ,v) \rangle\\
&=& K((i\alpha ,u),(j \beta ,v))=0
\end{eqnarray*}
since neither $i\alpha = j \beta \gamma$ nor $j \beta = i \alpha \gamma$ is possible.
As required for dilations we have
$\tilde{T}^*_i\lambda(0,u)=\lambda(0,T^*_iu)$ as may be seen as follows.
\begin{eqnarray*}
\langle \tilde{T}^*_i\lambda(0,u),\lambda(\beta,v)\rangle
&=&\langle \lambda(0,u),\tilde{T}_i\lambda(\beta,v)\rangle \\
&=& \langle \lambda(0,u),\lambda(i\beta,v)\rangle\\
&=& K((0,u),(i \beta,v))\\
&=& K((0,T^*_iu),(\beta,v))\\
&=& \langle \lambda(0,T^*_iu),\lambda(\beta,v)\rangle.
\end{eqnarray*}
Next we show that $\tilde{T}_i$ is a partial isometry by evaluating $\tilde{T}^*_i\tilde{T}_i \lambda(\alpha,u)$. By definition of $K$ we have
\[
\langle  \lambda(\alpha,u) , \tilde{T}^*_i\tilde{T}_i \lambda (\beta ,v) \rangle =
\left\{
\begin{array}{cc}
a_{i ,o(\alpha)} a_{i ,o(\beta)}\langle \lambda(\alpha,u) , \lambda (\beta ,v) \rangle &
\mbox{if ~}\alpha \neq 0 \text{ or } \beta \neq 0 \\
\langle u, Q_i v \rangle & \mbox{if ~} \alpha = \beta =0.
\end{array}
\right.
\]
Thus $\tilde{T}^*_i\tilde{T}_i \lambda(\beta,u) =a_{i ,o(\beta)}\lambda(\beta,u)$ if $\beta \neq 0$ and
\begin{eqnarray*}
\langle  \lambda(0 ,u), \tilde{T}^*_i \tilde{T}_i \lambda(0,v) \rangle
&=&
\langle u,Q_i v \rangle \\
&=&
\langle u , \sum_k a_{ik} T_k T_k^* v \rangle \\
&=&
\langle \lambda (0,u) , \sum_k a_{ik}  \lambda (k,T_k^*v) \rangle.
\end{eqnarray*}
Since $\langle \lambda (\alpha ,u) , \lambda (k , T_k^* v) \rangle = \delta_{k o(\alpha)} \langle u , (\uT^{\alpha})^* v \rangle$ for $\alpha \neq 0$ we have
\begin{eqnarray*}
\langle  \lambda(\alpha ,u), \tilde{T}^*_i \tilde{T}_i \lambda(0,v) \rangle
&=&
a_{i,o(\alpha)}\langle \lambda(\alpha ,u),\lambda(0,v) \rangle \\
&=&
\langle  \lambda(\alpha,u), \sum_k a_{ik} \lambda (k,T_k^* v) \rangle
\end{eqnarray*}
and therefore
$$
\tilde{T}^*_i \tilde{T}_i \lambda(\alpha,u)=\left\{
\begin{array}{cc}
\sum_k a_{ik} \lambda(k,T^*_k
u)&\mbox{if~}\alpha=0\\
a_{i o(\alpha)}\lambda(\alpha,u)& \mbox{otherwise.}
\end{array}
\right.
$$
It follows that $\tilde{T}_i \tilde{T}^*_i
\tilde{T}_i= \tilde{T}_i,$ i.e.\:the $\tilde{T}_i$
are partial isometries.
Finally minimality holds as
$$
\overline{\mbox{span}}\{\tilde{\uT}^\alpha \lambda (0,u):
 \alpha \in \tilde{\Lambda}, u \in \CH\}
=\overline{\mbox{span}}\{\lambda (\alpha,u): \alpha
\in \tilde{\Lambda}, u \in \CH\}=\tilde{\CH}.
$$

\medskip

\noindent
\textsc{Second proof using Popescu's method}

\medskip


Recall that $\uL$, $\uS$ denote the $n$-tuples of creation operators on
$\Gamma(\C^n)$, $\Gamma_A$ respectively.

\begin{Definition}{ \em
For a contractive tuple $\uT=(T_1,\cdots,T_n)$ on a Hilbert space $\CH$ the operator
$\Delta_{\uT}=(I-\sum^n_{i=1}T_iT^*_i)^{\frac{1}{2}}$ is called
{\em defect operator \/} of $\uT.$ If $\sum_{\alpha \in
\Lambda^m}\uT^\alpha(\uT^\alpha)^*$ converges to zero in the
strong operator topology as $m$ tends to infinity then this tuple
is said to be {\em pure.} }
\end{Definition}

Now let $\uT$ be a pure tuple on $\H$ satisfying $A$-relations. Similarly as in [Po1]
\begin{equation}
Kh= \sum_{\alpha\in \tilde{\Lambda}_A} e^\alpha \otimes \Delta_{\uT}(\uT^\alpha)^*h
\end{equation}
defines an isometry
$K:\CH\to
\Gamma_A \otimes \overline{\Delta_{\uT}(\CH)}$
such that $\uT^\alpha=K^*(\uS^\alpha \otimes
I)K$. Moreover for all $\alpha \in \tilde{\Lambda}, ~~~ S^*_i \otimes I$
leaves the range of $K$ invariant and
$$
\overline{\text{span}}\{(S_i \otimes I) Kh : i=1, \ldots ,n ,\; h \in \CH \} =\Gamma_A \otimes \overline{\Delta_{\uT}(\CH)}
$$
since $((I-\sum S_iS^*_i)\otimes I) Kh= \omega \otimes \Delta_{\uT}h$
and $\overline{\mbox{span}}\{\uS^\alpha \omega: \alpha \in
\tilde{\Lambda}_A \}= \Gamma_A$, the tuple
$(S_1 \otimes I,\cdots, S_n \otimes I)$ is the minimal Cuntz-Krieger dilation of $\uT$.
Note that for $\alpha, \beta \in \tilde{\Lambda}_A$ and $P_0 = I -\sum_i S_iS_i^*$ we have
\begin{eqnarray}
K^*[\uS^\alpha P_0(\uS^\beta)^* \otimes I]Kh
& = & K^*[\uS^\alpha P_0(\uS^\beta)^*\otimes
I] (\sum_\gamma \uS^\gamma \omega \otimes \Delta_{\uT}(\uT^\gamma)^*h) \nonumber \\
& = & K^* (\sum_\gamma \uS^\alpha P_0(\uS^\beta)^*\uS^\gamma \omega \otimes \Delta_{\uT}(\uT^\gamma)^*h)\\
& = & K^* ( \uS^\alpha \omega \otimes \Delta_{\uT}(\uT^\beta)^*h)=
\uT^\alpha \Delta_{\uT}^2 (\uT^\beta)^* h. \nonumber
\end{eqnarray}
As in [Po1] starting with a contractive tuple $\uT$ on a Hilbert space $\CH,$
the tuple $r\uT=(rT_1,\cdots,rT_n)$ is pure for $0 < r <1$. By (3.1) there
is an isometry $K_r:\CH \to \Gamma_A \otimes \overline {\Delta_r(\CH)}$ defined by
\begin{equation}
K_r h= \sum_\alpha e^\alpha \otimes
\Delta_r((r\uT)^\alpha)^*h,
\end{equation}
where $\Delta_r=(I-r^2 \sum T_iT^*_i)^{\frac{1}{2}}$. From this we obtain a unital
completely positive map $\psi_r:C^*(\uS) \to B(\H)$ defined by
$\psi_r(X)=K^*_r(X\otimes I)K_r, ~~X\in C^*(\uS)$. As the family of maps
$\psi_r$, where $0<r<1$, is uniformly bounded, $\psi_r$ converges pointwise.  Taking the limit as
$r$ increases to $1,$ we get a unique unital completely positive
map $\theta$ from $C^*(\uS)$ to $B(\CH)$
satisfying
\begin{equation}
\theta (\uS^\alpha (\uS^\beta)^*)=\uT^\alpha(\uT^\beta)^*
~~\mbox{for~} \alpha, \beta \in \tilde{\Lambda}_A.
\end{equation}
Once we have this map we can use a minimal Stinespring dilation
$\pi_1: C^* (\uS) \to B(\tilde{\CH})$ of $\theta$ such that
$$
\theta (X)=P_{\CH}\pi_1(X)|_{\CH}~~~ \forall X \in C^*(\uS)
$$
and $\overline{\mbox{span}}\{\pi_1(X)h: X \in C^*(\uS), h\in
\CH\}=\tilde{\CH}$. The tuple
$\tilde{\uT}=(\tilde{T}_1,\cdots,\tilde{T}_n)$ where
$\tilde{T}_i=\pi_1(S_i),$ is the minimal Cuntz-Krieger dilation
of $\uT$ which is unique up to unitary equivalence.
$\tilde{\uT}$ consists of partial isometries with
orthogonal ranges satisfying $A$-relations.

We remark that if $\uR$ is a tuple consisting of partial isometries
orthogonal ranges satisfying the condition
$$
R^*_iR_i=I-\sum^n_{j=1}(1-a_{ij})R_jR^*_j
$$
(e.g. a Cuntz-Krieger dilation) then the completely positive map $\Theta$ in (3.4)
mapping $\uS^{\alpha} (\uS^*)^{\beta}$ to $\uR^{\alpha} (\uR^*)^{\beta}$
is $*$-homomorphisms because the $S_i$'s and $R_i$'s have
orthogonal ranges and for $1 \leq i \leq n$
\begin{eqnarray*}
\Theta (S^*_iS_i)& = & \Theta (I-\sum_j (1-a_{ij})S_jS^*_j )
= I-\sum_j (1-a_{ij})R_jR^*_j\\
&=& R^*_iR_i = \Theta (S^*_i) \Theta(S_i).
\end{eqnarray*}

\medskip

\noindent
\textsc{The algebra generated by a Cuntz-Krieger dilation}

\medskip

The tuple obtained from the above constructions satisfies Cuntz-Krieger relations,
that is:
\begin{equation}
\tilde{T}^*_i\tilde{T}_i=I-\sum_j
(1-a_{ij})\tilde{T}_j\tilde{T}^*_j.
\end{equation}
We consider the $C^*$-algebra generated by such tuples.

First consider the $C^*$-algebra generated by left creation operators on $A$-Fock space.
Since for any $\alpha, \beta \in \tilde{\Lambda}_A$ the rank one operator
$\eta \mapsto \langle \uS^\beta \omega ,\eta \rangle \uS^\alpha \omega$ on $\Gamma_A$ can be
written as $\uS^\alpha(I-\sum S_iS^*_i)(\uS^\beta)^*=\uS^\alpha P_0 (\uS^\beta)^*$ and they
span the subalgebra of compact operators in $C^*(\uS),$ we
conclude that $C^*(\uS)$ also contains all compact operators.
For $\tilde{\uT}=\pi_1 (\uS)$ the Hilbert space
$\tilde{\CH}$ can be decomposed as $\tilde{\CH}=\tilde{\CH}_C\oplus
\tilde{\CH}_N$, where
$$
\tilde{\CH}_C:=\overline{\mbox{span}}\{\pi_1(X)h: h \in \tilde{\CH},
X \in C^*(\tilde{\uT}) \mbox{~and compact}\}
$$
is bi-invariant with respect to the $\tilde{T}_i$'s, that is,
invariant with respect to $\tilde{T}_i$ and $\tilde{T}^*_i$ for all $i$.
Thus $\pi_1$ can be decomposed as $\pi_{1C}\oplus \pi_{1N}$, where
$\pi_{1C}(X)=P_{\tilde{\CH}_C}\pi_1(X)P_{\tilde{\CH}_C}$ and
$\pi_{1N}(X)=P_{\tilde{\CH}_N}\pi_1(X)P_{\tilde{\CH}_N}.$ As
$\pi_{1N}$  annihilates  compacts,  $\pi_{1N}(I-\sum
S_iS^*_i)=\pi_{1N}(P_0)=0$, hence
$(\pi_{1N}(S_1),\cdots,\pi_{1N}(S_n))$ satisfy Cuntz-Krieger
relations (in particular $A$-relations) and  generate a Cuntz-Krieger algebra.

Let $\CK$ be the range of $\pi_1(P_0)$ then $\pi_{1C}(\uS^{\alpha})k \mapsto e^{\alpha} \otimes k$
extends to a unitary equivalence between $\pi_{1C}(\uS)$ and $\uS \otimes I$ on $\Gamma_A \otimes \CK$
so that $\CK$ is a wandering subspace for $\tilde{\uT}$ generating $\tilde{\CH}_C$.

The isometry in Stinespring's Theorem is of the form $V=\left[ \begin{array}{c} V_1\\
V_2 \end{array}\right]$ such that $V_1$ maps $\CH$ to
$\Gamma_A \otimes \CK$ and $V_2$ maps $\CH$ to $\tilde{\CH}_N$. Now
for $\alpha, \beta \in \tilde{\Lambda}$
\begin{eqnarray*}
\uT^\alpha \Delta_{\uT}^2 (\uT^\beta)^* h & = & \theta (\uS^\alpha
(I-\sum_i S_i S^*_i)(\uS^\beta)^*)(h)\\
& = &V_1^*[\uS^\alpha P_0 (\uS^\beta)^*\otimes I]V_1(h)
+  V_2^*[\pi_{1N}(\uS^\alpha P_0(\uS^\beta)^*)]V_2(h)\\
& = & V_1^*[\uS^\alpha P_0(\uS^\beta)^*\otimes I]V_1(h)
\end{eqnarray*}
as $\pi_{1N}$ annihilate compacts. Comparison with identity (3.2) shows that
$V_1$ may be taken to be $K$. Hence $\CK:=\overline{\Delta_{\uT} (\CH)}$.

In fact given just a tuple $\uR$ verifying
$$
R^*_i R_i=I-\sum^n_{j=1}(1-a_{ij})R_jR^*_j
$$
it is clear that we will always obtain
a decomposition of this type as the minimal Cuntz-Krieger dilation of such a
tuple is the tuple itself. Such a decomposition is called {\em Wold
decomposition.\/}

Using arguments similar to Theorem 1.3
in [Po1] we conclude that
$$
\tilde{\CH}_N=\bigcap^\infty_{m=0} \overline{\mbox{span}}\{\tilde{\uT}^\alpha h:
h \in \tilde{\CH}, |\alpha|=m \}.
$$

\begin{Corollary}
Suppose $\hat{\uT}$ is the minimal isometric dilation of a contractive tuple $\uT$ satisfying $A$-relations.
\begin{enumerate}
\item $\textup{rank~} (I -\sum_i \hat{T}_i \hat{T}^*_i)=\textup{rank} (I -\sum_i
\tilde{T}_i\tilde{T}^*_i)=\textup{rank} (I -\sum_i T_iT^*_i).$
\item $\lim_{k \to \infty} \sum_{|\alpha|=k} \tilde{\uT}^\alpha
(\tilde{\uT}^\alpha)^* = P_{\tilde{\CH}_N}.$
\end{enumerate}
\end{Corollary}
\noindent {\sc Proof:} Clear. \qed

\medskip

Now we would like to see how the minimal Cuntz-Krieger dilation
and minimal isometric dilation are related. The following is an
analogue of Theorem 13 in [BBD] for the maximal $A$-relation
piece.

\begin{Theorem}
Let $\uT$ be a contractive $n$-tuple on a Hilbert space $\CH$
satisfying $A$-relations. Then the maximal $A$-relation piece of
the standard noncommuting dilation $\hat{\uT}$ of $\uT$ is a realization of
the minimal Cuntz-Krieger dilation $\tilde{\uT}$ of $\uT$.
\end{Theorem}

\noindent {\sc Proof:}
Let $\theta:C^*(\uS)\to B(\CH)$ be the unital
completely positive map as in equation (3.4), $\pi_1$
the corresponding minimal Stinespring dilation and
$\tilde{T}_i=\pi_1(S_i)$ as before. Since the standard tuple $\uS$ on
$\Gamma_A$ is also a contractive tuple, there is a unique
completely positive map $\varphi$ from the $C^*$-algebra
$C^*(\uL)$ generated by the left creation operators to $C^*(\uS),$
satisfying
$$
\varphi(\uL^\alpha(\uL^\beta)^*)=\uS^\alpha(\uS^\beta)^*
\mbox{~~for ~~}\alpha,\beta \in \tilde{\Lambda}.
$$
Thus $\psi$
as defined before on $C^*(\uL)$ satisfies $\psi=\theta \circ
\varphi.$ Let the minimal Stinespring dilation of $\pi_1 \circ
\varphi$ be the $*$-homomorphism $\pi:C^*(\uL)\to B(\CH_1)$ for
some Hilbert space $\CH_1=\overline{\mbox{span}}\{\pi(X)h:X\in
C^*(\uL), h\in \tilde{\CH} \}$ such that
$$ \pi_1 \circ \varphi(X)= P_{\tilde{\CH}} \pi
(X)|_{\tilde{\CH}}~~~\forall X \in C^*(\uL).$$ In the following
commuting diagram

\hskip1.5in
\begin{picture}(200,115)

\put(0,20){$C^*(\uL)$}

\put(40,21){$\longrightarrow$}

\put(70,20){$C^*(\uS)$}

\put(110,21){$\longrightarrow$}

\put(140,20){$\mathcal{B}(\CH)$ \hspace{0.4cm} $\uT$}

\put(140,60){$\mathcal{B}(\tilde{\CH})$ \hspace{0.4cm} $\underline{\tilde{T}}$}

\put(140,100){$\mathcal{B} (\CH_1)$ \hspace{0.4cm} $\hat{\uL}=\hat{\tilde{\uT}}$}

\put(40,35){\vector(4,3){83}}

\put(110,33){\vector(4,3){23}}

\put(150,40){$\downarrow $}

\put(150,80){$\downarrow $}

\put(45,10){$\varphi$}

\put(115,10){$\theta$}

\put(110,46){$\pi_1$}

\put(70,72){$\pi$}

\end{picture}

\noindent all horizontal arrows are unital completely positive
maps, down arrows are compressions and diagonal arrows are minimal
Stinespring dilations. Let $\hat{L}_i=\pi(L_i)$ and
$\hat{\uL}=(\hat{L}_1,\cdots,\hat{L}_n)= \hat{\tilde{\uT}}$. We will first
show that $\tilde{\uT}$ is the maximal A-relation piece of
$\hat{\uL}$ and then show that $\hat{\uL}$ is the standard
noncommuting dilation of $\uT.$

To this end we use the presentation of the minimal isometric dilation $\hat{\uL}$ which was given
by Popescu [Po1]. In the one variable case it was given by Sch\"affer c.f.\;also [BBD]. Define
$D:\underbrace{\tilde{\CH} \oplus \cdots \oplus
\tilde{\CH}}_{n-\mbox{copies}}
\to
\underbrace{\tilde{\CH} \oplus \cdots \oplus \tilde{\CH}}_{n-\mbox{copies}}$
by
$$
D^2=[\delta_{ij}I-\tilde{T}^*_i\tilde{T}_j]_{n \times
n}=[\delta_{ij}(I-\tilde{T}^*_i\tilde{T}_i)]_{n
\times n}
$$
as used by Popescu. Note that $D^2$ is a projection as all $\tilde{T}_i$'s
are partial isometries hence $D^2=D.$ Let
$\CD$ denote the range of $D.$ We identify $\underbrace{\tilde{\CH}
\oplus \cdots \oplus \tilde{\CH}}_{n-\mbox{copies}}$ with $\C^n\otimes
\tilde{\CH}$
and hence $(h_1,\cdots, h_n)$ with
$\sum^n_{i=1}e_i\otimes h_i$ and $\C \omega
\otimes \CD$ with $\CD.$
$$
D(h_1,\cdots, h_n)=D(\sum^n_{i=1}e_i\otimes h_i)=\sum^n_{i=1}
e_i\otimes (I-\tilde{T}^*_i\tilde{T}_i)h_i.
$$
For $h\in \tilde{\CH}, d_\alpha
\in \CD, 1 \leq i \leq n,$ the standard noncommuting dilation
$\hat{\uL}=(\hat{L}_1,\cdots,\hat{L}_n)$ is given by
$$
\hat{L}_i(h\oplus \sum_{\alpha \in \tilde{\Lambda}} e^\alpha \otimes
d_{\alpha}) =\tilde{T}_i h \oplus D(e_i\otimes h)\oplus e_i
\otimes(\sum_{\alpha \in \tilde{\Lambda}} e^\alpha \otimes
d_\alpha)
$$
on the dilation space $\CH_1=\tilde{\CH} \oplus
(\Gamma(\C^n)\otimes \CD)$. As $\tilde{\uT}$ satisfies
$A$-relations and $\hat{L}^*_i$ leaves $\tilde{\CH}$ invariant,
$\tilde{\CH} \subseteq (\CH_1)_A(\hat{\uL})$. To show the reverse inclusion suppose that there exists a non-zero
$z \in \tilde{\CH}^ \bot \cap (\CH_1)_A(\hat {\uL})$. $z$ can be written as $0\oplus
\sum_{\alpha \in \tilde{\Lambda}} e^\alpha \otimes z_\alpha$ such
that $z_\alpha \in \CD$. Since $\langle \omega \otimes z_\alpha,(\hat{\uL}^\alpha)
^*z\rangle=\langle e^\alpha\otimes z_\alpha, z\rangle=\langle
z_\alpha,z_\alpha\rangle$ and $(\hat{\uL}^\alpha) ^*z \in
(\CH_1)_A(\hat {\uL})$, we can assume $\| z_0\|=1$ without loss of
generality. Also $z_0=D(h_1,\cdots,h_n)$ for some $h_i \in
\tilde{\CH}$ as projections have closed ranges.
Now consider
\begin{eqnarray*}
\sum^n_{i,j=1}(\hat{L}_i\hat{L}_j-a_{ij}\hat{L}_i\hat{L}_j)\tilde{T}^*_jh_i
&=&\sum^n_{i,j=1}(\tilde{T}_i\tilde{T}_j-a_{ij}\tilde{T}_i\tilde{T}_j)\tilde{T}
^*_jh_i+\sum^n_{i=1}D(e_i
\otimes\sum^n_{j=1}(1-a_{ij})\tilde{T}_j\tilde{T}^*_jh_i)\\
& &+ \sum^n_{i,j=1}(1-a_{ij})e_i\otimes D(e_j \otimes
\tilde{T}^*_jh_i)\\
&=&0+\sum^n_{i=1}D[e_i\otimes(I-\tilde{T}^*_i\tilde{T}_i)h_i]+x\\
&=&D^2(h_1,\cdots, h_n) + x = \tilde{z}_0+x.
\end{eqnarray*}
where
$x=\sum^n_{i,j=1}(1-a_{ij})e_i\otimes D(e_j\otimes
\tilde{T}^*_jh_i)$.
Thus $\langle z,\tilde{z}_0+x \rangle=0$ by Lemma 4. Moreover
\begin{eqnarray*}
\|x\|^2&=&\| \sum^n_{i,j=1}(1-a_{ij})e_i\otimes
 D(e_j \otimes \tilde{T}^*_jh_i)\|^2\\
&=&\sum^n_{i=1}\langle \sum^n_{j=1}(1-a_{ij})D(e_j\otimes
\tilde{T}^*_jh_i),
\sum^n_{j'=1}(1-a_{ij'})e_{j'}\otimes \tilde{T}^*_{j'}h_i\rangle\\
&=&\sum^n_{i=1}\langle
\sum^n_{j=1}(1-a_{ij})(I-\tilde{T}^*_j\tilde{T}_j)\tilde{T}^*_jh_i,
\sum^n_{j=1}(1-a_{ij})\tilde{T}^*_jh_i\rangle\\
&=&\sum^n_{i=1}\langle
\sum^n_{j=1}(1-a_{ij})(\tilde{T}^*_j-\tilde{T}^*_j)h_i,\sum^n_{j=1}
(1-a_{ij})\tilde{T}^*_jh_i\rangle=0,
\end{eqnarray*}
i.e. $x=0$ and therefore $\|\tilde{z}_0\|^2=\langle z,\tilde{z}_0\rangle=0$
which is a contradiction. Hence $z=0$ which implies
$\CH=(\CH_1)_A(\hat{\uL})$.\\

To show finally that $\hat{\uL}$ is standard note that
$$
\CH_1=\overline{\text{span}}\{ \hat{\uL}^{\alpha} x : \alpha \in \Tilde{\Lambda}_A , \; x \in \tilde{\CH} \}
$$
and
$$
\tilde{\CH}=\overline{\text{span}}\{ \tilde{\uT}^{\alpha} z : \alpha \in \Tilde{\Lambda}_A , \; z \in \CH \},
$$
moreover $P_{\tilde{\CH}} \hat{\uL}^{\alpha} =\tilde{\uT}^{\alpha}$ and $P_{\CH} \hat{\uL}^{\alpha} =\uT^{\alpha}$
for $\alpha \in \tilde{\Lambda}_A$ by assumption. Hence
\begin{eqnarray*}
\CH_1
&=&
\overline{\text{span}}\{ \hat{\uL}^{\alpha} x : \alpha \in \Tilde{\Lambda}_A , \; x \in \tilde{\CH} \} \\
&=&
\overline{\text{span}}\{ \hat{\uL}^{\alpha} \tilde{\uT}^{\beta} z : \alpha , \beta \in \Tilde{\Lambda}_A , \; z \in \CH \} \\
&=&
\overline{\text{span}}\{ \hat{\uL}^{\alpha} P_{\tilde{\CH}} \hat{\uL}^{\beta}  z : \alpha , \beta \in \Tilde{\Lambda}_A , \; z \in \CH \} \\
& \subseteq &
\overline{\text{span}}\{ \hat{\uL}^{\alpha} \hat{\uL}^{\beta} z : \alpha , \beta \in \Tilde{\Lambda}_A , \; z \in \CH \}= \CH_1.
\end{eqnarray*}
\qed

\medskip

In the same way one can show that similar results holds even for $q$-commuting tuples (c.f.\:[BBD]).
To keep the presentation simpler we have worked with the above
special case. The following example illustrates the forgoing results.

\begin{Example} {\em
 For $\CH=\C^2,$ let $T_1= \left(
\begin{array}{cc} 0 &1 \\
0 & 0
\end{array}\right)$ and $T_2= \left( \begin{array}{cc} 0 &0 \\ 1 & 0
\end{array}\right)$.  Then one observes that $\uT$ satisfies
$A$-relations for the matrix  $A= \left(
\begin{array}{cc} 0 &1 \\
1 & 0
\end{array}\right)$, $T_1T^*_1+T_2T^*_2=I$ and the $T_i$'s are
partial isometries. Further $D$ used in the above Theorem turns out to
be
$$D=\left( \begin{array}{cccc} 1 & & & \\  & 0
& & \\ & & 0 & \\ &&& 1 \end{array}\right).$$ Let us denote  the
two basis vectors in the range of $D$ corresponding to the entries 1
appearing in $D$ by $f_1$ and $f_2$ such that
$$D(e_1 \otimes (a_1, a_2) + e_2 \otimes (b_1, b_2))= a_1 f_1 +b_2 f_2$$
for all $a_1,a_2, b_1, b_2 \in \C.$ The dilation space for the
minimal isometric dilation $\uV=(\tilde{V}_1,\tilde{V}_2)$ of $\uT$ is
$\CH \oplus \Gamma (\C^n) \otimes
\CD $ where $\CD$ is the range of $D.$
$$
\tilde{V}_1\tilde{V}_1(a_1,a_2)=(0,0)+\omega \otimes(a_2,0) +e_1 \otimes
(a_1,0),
$$
$$
\tilde{V}_2\tilde{V}_2(a_1,a_2)=(0,0)+\omega \otimes(0,a_1) +e_2 \otimes
(0,a_2).
$$ As $a_1, a_2$ are arbitrary using the above equations
together with the description of the maximal $A$-relation piece
from Lemma 4 we get $\tilde{\CH}=\CH.$ }
\end{Example}

\end{section}

\begin{section}{Representations of Cuntz-Krieger Algebras}
\setcounter{equation}{0}

In general a Cuntz-Krieger algebra $\CO_A$ admits many inequivalent representations.
When $\uT=(T_1,\cdots, T_n)$ is a tuple on the Hilbert space $\CH$
satisfying $A$-relations and $\sum^n_{i=1}T_iT^*_i=I,$ the minimal
Cuntz-Krieger dilation
$\tilde{\uT}=(\tilde{T}_1,\cdots,\tilde{T}_n)$ is such that
$C^*(\tilde{\uT})$ is a Cuntz-Krieger algebra.
If $\uT$ is moreover commuting then the unital
completely positive map $\rho_{\uT}:\CO_A \to
C^*(\tilde{\uT})$ given by $\rho_{\uT}(s_i)=\tilde{T}_i$ is a
representation of $\CO_A.$ We will classify such representations
here.

For a tuple $\uR=(R_1,\cdots,R_n)$ on a Hilbert space $\CK,$ we
use the concept of  {\em maximal commuting piece} and the
space $\CK^c(\uR)$ as defined before Lemma 6 and section 2 of [BBD]. We refer to
$\K^c(\uR)$ as {\em maximal commuting subspace.}
\begin{Definition} {\em
\begin{enumerate}
\item A commuting tuple $\uT=(T_1,\cdots,T_n)$ is called {\em spherical
unitary \/} if $\sum T_iT^*_i=I$ and all $T_i$'s are normal.
\item A representation $\rho$ of $\CO_A$ on $B(\CK)$ for some
Hilbert space $\CK,$ is said to be {\em spherical \/} if
$R_i=\rho(s_i), 1\leq i \leq n$ and $\CK=\{ \uR^\alpha k: k \in
\CK^c(\uR) \mbox{~and~}  \alpha \in \tilde{\Lambda} \}.$
\end{enumerate} }
\end{Definition}

If $\uT$ is a spherical unitary then by Fuglede's Theorem $C^*(\uT)$ is
a commutative $C^*$-algebra i.e.\;$T_i^*$ and $T_j$ also commute for all $i$ and $j$.

\begin{Definition} {\em The {\em maximal commuting $A$-subspace \/} of a
$n$-tuple of isometries $\uV$ with orthogonal ranges is
defined as the intersection of its maximal commuting subspace and
maximal $A$-relation subspace. The $n$-tuple obtained by
compressing each $V_i$ to the maximal commuting $A$-subspace
is called {\em maximal commuting $A$-piece.}}
\end{Definition}

\begin{Remark} {\em Making use of Lemma 4, it follows that
$$
\CK_A \cap \CK^c= (\CK_A)^c = (\CK^c)_A
$$
i.e.\;the maximal commuting $A$-subspace of a $n$-tuple
 is in fact the maximal commuting subspace of
 the maximal $A$-relation piece or
 the  maximal $A$-relation subspace of the
 maximal commuting piece. }
\end{Remark}

Let $P_0=1$ on $\C$ and $P_m$ be the projection $\frac{1}{m!}\sum_{\sigma \in \CS_m} U^m_\sigma$
acting on $(\C^n)^{\otimes^m}$ where
$$
U^m_\sigma (y_1\otimes \cdots \otimes y_m)=y_{\sigma^{-1}
(1)}\otimes \cdots \otimes y_{\sigma^{-1}(m)},
$$
$y_i\in \C^n.$ We denote $\oplus^\infty_{m=0}P_m$ by $P^s.$
Given $A$-relations define
$\Lambda^m_{sA} = \{ \alpha \in \tilde{\Lambda}:  \mbox{~either~} |\alpha|=m>1 \mbox{~and~}
a_{\alpha_i\alpha_j}=1 \mbox{~for~} 1 \leq i \neq j\leq m,
\mbox{~or~}
|\alpha|\leq 1\} \subseteq \Lambda^m_A$.

\begin{Definition} {\em The subspace of $\Gamma_A$ defined by
$$\overline{\mbox{span}}\{ P^s e^\alpha: \alpha \in \Lambda_{sA}^m \}.$$ is
called {\em
commuting $A$-Fock space \/} and denoted by $\Gamma_{sA}.$ }
\end{Definition}

To see that $\Gamma_{sA}$ is the maximal commuting $A$-subspace of
$\uL$ we first note that the maximal commuting $A$-subspace of $\uL$ is the
intersection of symmetric Fock space $\Gamma_s (\C^n)$ (c.f.\;[BBD]) and
the maximal $A$-relation subspace of $\uL.$ Also
$$
\Gamma_s (\C^n)=\overline{\mbox{span}}\{ P^s e^\alpha : \alpha \in
\tilde {\Lambda}\}.
$$
Suppose $\alpha \in \Lambda^m$ and for all $1 \leq k \neq l \leq m,
a_{\alpha_k \alpha_l}=1$ then for $h \in \Gamma (\C^n)$ and all $i,j$
$$
\langle P^s e^\alpha, \uL^\beta (L_i L_j - a_{ij} L_i L_j)h \rangle=0.
$$
So, from the definition it is clear that
$$
\Gamma_{sA} \subseteq \Gamma_s(\C^n) \cap \Gamma_A.
$$

Let $\hat{P}$ denote the projection onto $\Gamma_s(\C^n) \cap \Gamma_A$
and let $z \in \Gamma_s(\C^n) \cap \Gamma_A$
be arbitrary. Suppose
$\alpha \in \Lambda^m$ is such that $\langle e^\alpha,z \rangle$ is not equal to $0.$
As $z \in \Gamma_A$, it follows that $\alpha \in \Lambda_A^m$.  Further
for any $\sigma \in S_m$
\begin{eqnarray*}
\langle U^m_\sigma e^\alpha, z\rangle &=& \langle U^m_\sigma e^\alpha, \hat{P} z\rangle =
\langle \hat{P} U^m_\sigma e^\alpha,  z\rangle\\
&=& \langle \hat{P} e^\alpha, z\rangle = \langle e^\alpha, z\rangle.
\end{eqnarray*}
Thus $ \langle U^m_\sigma e^\alpha, z\rangle$ is not equal to $0.$ This implies that
$\alpha \in \Lambda_{sA}^m$ and hence $z \in \Gamma_{sA}.$ We conclude that
$\Gamma_{sA}$ is the maximal commuting $A$-subspace of
$\uL.$

Also we would like to remark that the Fermionic Fock space (c.f.\:[De])
$\Gamma_a(\C^n)$ is the intersection of the maximal $q$-commuting
subspace (defined in [De]) and the maximal $A$-relation subspace with respect to the
following $q=(q_{ij})_{n \times n}$ and $A=(a_{ij})_{n \times n}:$

$$
q_{ij}=\left\{ \begin{array}{cc}
1&\mbox{if~} i = j\\
-1& \mbox{ otherwise}\end{array} \right.
\mbox{~and~}a_{ij}=\left\{ \begin{array}{cc}
0& \mbox{if~} i =j\\
1& \mbox{otherwise.} \end{array}\right.
$$
It is easy to see this using arguments similar to that we use to show that $\Gamma_{sA}$ is
the maximal commuting $A$-subspace with respect to $\uL.$
In other words, the Fermionic Fock space
$\Gamma_a(\C^n)$ is the maximal piece for the set of polynomials:
$$
p_{1ij}(\uz)= z_j z_i -q_{ij}z_i z_j\mbox{~and~} p_{2ij}(\uz)=z_i z_j
- a_{ij} z_i z_j \forall 1\leq i,j \leq n.
$$

Notice that $L^*_i$ leaves $\Gamma_{sA}$ invariant as $S^*_i$ leaves $\Gamma_{sA}$ invariant. Let
the compression of $L_i$ onto $\Gamma_{sA}$ be denoted by $W_i$, i.e.\;$\uW$
is the maximal commuting $A$-relation piece of $\uL$.
Suppose $\alpha \in \Lambda_{sA}^m,$ and when $ |\alpha| > 0$ let
$\alpha=(\alpha_1, \cdots, \alpha_m)$ where $m=|\alpha|.$ The operator
$W_i$ turns out to be
$$
W_iP^se^\alpha=\left\{ \begin{array}{cc} e_i & \mbox{~if~} |\alpha|=0 \\
P^se_i \otimes e^\alpha &  \mbox{~if~} a_{i \alpha_j}a_{\alpha_j i}=1,
\forall 1\leq j \leq m \\
0 & \mbox{~otherwise.} \end{array} \right.
$$
Form this it follows that
$\uW=(W_1,\cdots,W_n)$ is the maximal commuting piece satisfying $A$-relations of $\uL.$
Let us denote the maximal commuting piece of $\uL$ on $\Gamma (\C^n)$
by $\uC=(C_1, \cdots,C_n)$. (These are just the creation operators on symmetric Fock space.)
Then for $\alpha \in \Lambda_{sA}^m, \alpha=(\alpha_1,\cdots, \alpha_m),m>1$
the commutators verify
$$
[W_i,W^*_i]P^se^\alpha=\left\{ \begin{array}{cc} [C_i,C^*_i]P^se^\alpha
 & \mbox{~if~} a_{i \alpha_j}a_{\alpha_j i}=1, \forall 1\leq j \leq m \mbox{~or if~}\alpha=0 \\
\frac{1}{m!}P^se^\alpha & \mbox{if $\alpha_j=i$ for some $1\leq j \leq m$
and $a_{ii}=0$}\\
0 &\mbox{otherwise.}\end{array}\right.
$$
It is known that $[C_i,C^*_i]$ is compact for all $i$ (c.f.\:[Ar3, 5.3]), so
the  $[W_i,W^*_i]$'s are compact. Clearly, the vacuum vector is
contained in $\Gamma_{sA}$ and $I-\sum W_iW^*_i$ is the projection
on to the vacuum space. It contains all the rank one operators of the type
$\mu \to \langle \uW^{\alpha} \omega, \mu \rangle \uW^\beta \omega$ on $\Gamma_{sA}$
as those can be written as $\uW^\alpha ( I-\sum W_iW^*_i) (\uW^{\beta})^*$.
As these rank one operators span the subalgebra of compact operators, we conclude
that $C^*(\uW)$ contains the subalgebra of all compacts $\CK$ as an ideal. Since the image of
$\uW$ in $C^*(\uW)/ \CK$ is a spherical unitary it follows from Fuglede's Theorem that also $[W_i,W_j^*]$
where $i \neq j$ must be compact and we also conclude that
$$
C^*(\uW)=\overline{\mbox{span}}\{ \uW^\alpha (\uW^\beta)^*: \alpha,
\beta \in \tilde{\Lambda} \}.
$$

For a commuting pure tuple $\uT$ satisfying $A$-relations, with easy
computations it can be seen that the range of the isometry $K_r: \CH \to
\Gamma_A\otimes \overline{\Delta_{\uT}(\CH)}, 1 \leq r \leq 1,$ defined
in equation (3.4) is contained in $\Gamma_{sA}\otimes
\overline{\Delta_{\uT}\CH}$ and we obtain a unital completely positive map
$\phi: C^*(\uW) \to B(\CH)$ defined as strong operator topology limit of
$K^*_r(. \otimes I)K_r$ as $r$ increases to $1.$
 Let $\pi_0: C^*(\uW) \to B(\CH_0)$
be the minimal Stinespring dilation of $\phi$ for some Hilbert space $\CH_0$ and
$\check{W}_i=\pi_0(W_i)$ where $\CH_0=\overline{\mbox{span}}\{
\check{\uW}^\alpha h: \alpha \in \tilde{\Lambda}, h \in \CH \}.$

\begin{Definition}
{\em The above defined tuple $\check{\uW}=(\check{W}_1,\cdots, \check{W}_n)$ is said to be the {\em
standard commuting $A$-dilation of $\uT$.}}
\end{Definition}
\begin{Remark}
{\em It follows from Theorem 15 in [BBD] that for spherical
unitaries $\uT$ satisfying $A$-relation the maximal commuting piece
of the standard noncommuting dilation is $\uT.$ As $\uT$ satisfies
$A$-relations, it is clear that $\uT$  is also the maximal commuting
$A$-piece.}
\end{Remark}

So far for a commuting contractive tuple satisfying $A$-relations we have four types of standard minimal dilations: the isometric dilation, the Cuntz-Krieger dilation (or
$A$-dilation),
the commuting dilation and the commuting $A$-dilation. These are obtained by
considering Stinespring dilations of suitable completely positive maps on
$C^*(\uL), C^*(\uS), C^*(\uC),$ and $ C^*(\uW)$ respectively.
The last dilation is in a certain sense the intersection of the previous two.
The next Lemma, which is a generalization of Theorem 13 in [BBD] makes this
statement rigorous. This will be crucial for classifying certain types of
representations of Cuntz-Krieger algebras.
\begin{Lemma}
The maximal commuting piece of the minimal Cuntz-Krieger dilation
of a commuting tuple $\uT$ satisfying $A$-relations is the
standard commuting $A$-dilation.
\end{Lemma}
\noindent {\sc Proof:} Let the unital completely positive map
$\phi: C^*(\uW) \to B(\CH),$ $\pi_0$ and $\CH_0$ be as above. We
denote the operator $\pi_0(W_i)$ by $\check{W}_i$ and denote the
$n$-tuple $(\check{W}_1,\cdots,\check{W}_n)$ by $\check{\uW}.$
As $\uW$ is a contractive tuple satisfying $A$-relation,
there is a unital completely positive map $\eta: C^*(\uS) \to
C^*(\uW)$ such that $\eta (\uS^\alpha(\uS^
\beta)^*)=\uW^\alpha(\uW^\beta)^*.$ The completely positive map
$\theta$ as in equation (3.4) is equal to $\phi \circ \eta.$ Let
 $\tilde{\pi}_1$ be the minimal Stinespring dilation of $\pi_0
\circ \eta$ and $V_i=\tilde{\pi}_1 (S_i).$ We have the following commuting diagram.

\hskip1.5in
\begin{picture}(200,115)

\put(0,20){$C^*(\uS)$}

\put(40,21){$\longrightarrow$}

\put(70,20){$C^*(\uW)$}

\put(110,21){$\longrightarrow$}

\put(140,20){$\mathcal{B}(\CH)$}

\put(140,60){$\mathcal{B}(\CH_0)$}

\put(140,100){$\mathcal{B} (\tilde{\CH}_1)$}

\put(40,33){\vector(4,3){84}}

\put(110,32){\vector(4,3){23}}

\put(150,40){$\downarrow $}

\put(150,80){$\downarrow $}

\put(45,10){$\eta$}

\put(115,10){$\phi$}

\put(110,46){$\pi_0$}

\put(70,72){$\tilde{\pi}_1$}

\end{picture}\\
As before horizontal arrows are completely positive maps, diagonal
arrows are $*$-homomor-phism and down arrows are compressions.

Since $C^*(\uW)$  contains all compact
operators,  $\CH_0$ can again be decomposed as $\CH_{0C}\oplus
\CH_{0N}$ where $\CH_{0C}=\overline{\mbox{span}}\{ \pi_0(X)h:h \in
\CH, X\in C^*(\uW), X \mbox{~compact~} \}$ and $\CH_{0N}=\CH_0 \ominus \CH_{0C}$.
Correspondingly,
\[
\pi _0(X)=
\begin{pmatrix}
\pi_{0C}(X)&            \\
           & \pi_{0N}(X)
\end{pmatrix} ,
\]
where $\pi_{0C}(X)$ and $\pi_{0N}(X)$ are compressions of $\pi_0 (X)$ to
$\CH_{0C}$ and $\CH_{0N}$ respectively. Furthermore
$\CH_{0C}=\Gamma_{sA}\otimes \overline{\Delta_{\uT}(\CH)}$ and
$\pi_{0C}(X)=X\otimes I.$  Let $E_i=\pi_{0N}(W_i)$ and
$\uE=(E_1,\cdots,E_n).$ As $[W_i,W^*_i]$ and $I-\sum W_iW^*_i$ are
compacts, clearly $\uE$ consists of pairwise commuting normal operators i.e.
${\uE}$ is a spherical unitary satisfying A-relations.

From the properties of Popescu's Poisson transform
and $\Gamma_{sA},$ it follows that  $(W_1 \otimes
I,\cdots,W_n \otimes I)$ is the maximal commuting $A$-piece of its
standard noncommuting dilation  $(L_1 \otimes I,\cdots,L_n \otimes I).$
Also from Remark 21 we conclude that $\uE$ is the maximal commuting $A$-piece of
its standard noncommuting dilation. Using Remark 18 and Theorem 14
we observe that each of them is the maximal commuting piece of their
minimal Cuntz-Krieger dilation. Hence by Lemma 6, $\check{\uW}$ is the
maximal commuting piece of $\uV$.  From this using arguments
similar to Theorem 13 in [BBD] it can be shown that $\uV$
is the minimal Cuntz-Krieger dilation of $\check{\uW}.$ Hence the
Lemma follows. \qed

\medskip

 If a commuting contractive tuple $\uT$ also satisfies $A$-relations for
$A=(a_{ij})_{n \times n},$ then without loss of generality we may assume
$A$ to be symmetric, i.e., $A=A^t$. In this case $A$
is the adjacency matrix of a (nondirected) graph $G$ with
vertex set $\{ 1,2,\cdots,n\}$ and  set of edges $E=\{
(i,j):a_{ij}=1, 1 \leq i < j \leq n \}$. A vertex $i$ is said to
be a {\em zero vertex} if $a_{ii}=0$. Let us associate to
this graph a subset $M$ of $\{(z_1,\cdots,z_n):
\sum^n_{i=1}|z_i|^2=1\}$ defined as the set of elements satisfying $A$-relations, that is
$$
M=\{(z_1,\cdots,z_n): \sum^n_{i=1}|z_i|^2=1, z_iz_j=a_{ij}z_iz_j, 1 \leq i, j \leq n\}.
$$

The set $M$ can be described in the following way:
For a zero vertex $i,$ the
corresponding $z_i$ of any element of $M$ will always taken to be
zero.  For any element $(z_1,\cdots,z_n)$ of $M,$ some elements $z_{i_1},\cdots,z_{i_k}$ for different $1 \leq i_k \leq n$
can be simultaneously choosen to be non-zero if and only if $i_1,\cdots,i_k$ are nonzero vertices and form
vertices of an induced subgraph of $G$ which is also complete.

Let  $C^M_n$ be the $C^*$-algebra of continuous complex valued
functions on $M.$ Consider the tuple
$\uz=(z_1,\cdots,z_n)$ of co-ordinate functions $z_i$ in $C^M_n.$ To
any spherical unitary $\uR=(R_1,\cdots,R_n)$ satisfying $A$-relations there corresponds a
unique representation of $C^M_n$ mapping $z_i$ to $R_i.$ As for any
commuting tuple $\uT$ satisfying $A$-relations with $\sum T_iT^*_i=I,$ the standard
commuting dilation $\tilde{\uT}=(\tilde{T}_1,\cdots,\tilde{T}_n)$
is a spherical unitary (c.f.\;section 3 of [BBD]), we have a
representation $\eta_{\uT}$ of $C^M_n$ such that
$\eta_{\uT}(z_i)=\tilde{T}_i.$ From Theorem 14, it is easy to see
that if $\uD$ and $\uE$ are two commuting $n$-tuples of operators
satisfying the same $A$-relations (on not necessarily the same Hilbert
space), the corresponding representations $\rho_{\uD}$ and
$\rho_{\uE}$ of $\CO_A$ are unitarily equivalent if and only if
the representations $\eta_{\uD}$ and $\eta_{\uE}$ of $C^M_n$ are
unitarily equivalent.

Any $z=(z_1,\cdots,z_n)\in M$ satisfying $A$-relations as operator
tuple on $\C$ and is a spherical unitary. We can obtain a one
dimensional representation $\eta_z$ of $C^M_n$ which maps $f$ to
$f(z).$ Let $(V^z_1,\cdots, V^z_n)$ and $(S^z_1,\cdots,S^z_n)$ be
the standard noncommuting dilation and the minimal Cuntz-Krieger
dilation respectively of this operator tuple $z=(z_1,\cdots,z_n).$
The dilation space of the standard noncommuting dilation is
$$
\CH^z=\C\oplus (\Gamma(\C^n)\otimes \C^n_z),
$$
where $\C^n_z$ is the $(n-1)$-dimensional subspace of $\C^n$
orthogonal to $(\oz_1,\cdots,\oz_n)$ and
$$
V^z_i(h \oplus \sum_{\alpha} e^{\alpha} \otimes d_\alpha)=
a_i \oplus D(e_i \otimes h) \oplus e_i \otimes (\sum_\alpha
e^{\alpha} \otimes d_\alpha).
$$
Using the minimal Cuntz-Krieger
dilation $\uS^z$ we get a representation $\vartheta:\CO_A \to
C^*(\uS^z)$ mapping $s_i$ to $S^z_i.$ This is the GNS
representation of the Cuntz-Krieger state
$$
\us^\alpha (\us^\beta)^* \to \uz^\alpha \overline{(\uz^\beta)}.
$$
which exists by (3.4).
We call such states {\em Cuntz-Krieger states.}

\begin{Theorem}
Any spherical representation of $\CO_A$ (on a separable Hilbert space) can be written as direct
integral of GNS representations of Cuntz-Krieger states.
\end{Theorem}
\noindent{\sc Proof:} An arbitrary representation of $C^M_n$ is a
countable direct sum of multiplicity free representations. Also
any multiplicity free representation of $C^M_n$ can be seen as a map
which sends $g \in C^M_n$ to an operator which acts as multiplication
by $g$ on $L^2(M,\mu)$  for some finite Borel measure
$\mu$ on $M$ and the associated representation $\vartheta$ of $\CO_A$ can be expressed
as direct integral of representations $\vartheta_z$ with respect
to the measure $\mu$ acting on $\int \! \! \!\! \!  \oplus\CH^z\mu
(dz).$ Thus the Theorem follows.  \qed

\begin{Example}{\em Let $A=\left( \begin{array}{cccc}
 1 & 1& 1& 0\\
1 & 1& 0 & 0\\
1 & 1 & 0 & 0\\
1 & 0 & 0 & 1 \end{array} \right).$ Then any commuting contractive
$A$-relation tuple also satisfies $A'$-relations, where $A'$ is the symmetric matrix $\left(
\begin{array}{cccc}
 1 & 1& 1& 0\\
1 & 1& 0 & 0\\
1 & 0 & 0 & 0\\
0 & 0 & 0 & 1 \end{array} \right).$ Furthermore the set of vertices of
the corresponding graph is $\{ 1,2,3,4\},$ the set of edges is $E=\{(1,2),(1,3)
\}$ and $3$ is a zero vertex. Hence $M=[(\C^2\times
\{0\}^2)\cup(\{0\}^3 \times \C)]\cap \partial B_n$ where $\partial B_n=\{(z_1,\cdots,z_n):
\sum^n_{i=1}|z_i|^2=1\}. $}

\end{Example}

\begin{Corollary}
Any representation of $\CO_A$ can be decomposed as $\pi_s \oplus
\pi_t$ where $\pi_s$ is spherical representation and
$(\pi_t(s_1),\cdots,\pi_t(s_n))$ has trivial maximal commuting
piece.
\end{Corollary}
\noindent{\sc Proof:} Similar to proof of Theorem 19 in [BBD].\qed

\medskip

It also follows that for irreducible representations of $\CO_A,$ the maximal commuting
piece of $(\pi(s_1),\cdots,\pi(s_n))$ is either one dimensional or
trivial.
\end{section}

\begin{section}{Universal Properties and WOT-closed Algebras Related to
Minimal Cuntz-Krieger
Dilation}
\setcounter{equation}{0}

Assume $\tilde{\uT}$ to be the minimal Cuntz-Krieger dilation of a
contractive tuple $\uT$ satisfying $A$-relations.
 Define $C^*(\tilde{\uT})$ to be the unital $C^*$-algebra generated by $\tilde{\uT}.$ Clearly the linear map from $C^*(\tilde{\uT})$ to
$B(\CH)$ such that $ \tilde{\uT}^\alpha (\tilde{\uT}^\beta)^* \mapsto P_{\CH} \tilde{\uT}^\alpha (\tilde{\uT}^\beta)^* |_{\CH}=
\uT^\alpha (\uT^\beta)^*$ is a unital completely positive map. We will investigate some universal properties of minimal Cuntz-Krieger
dilations using methods employed by Popescu [Po5] for minimal isometric dilations.
Proposition 26
is a nonspatial characterization of the minimal Cuntz-Krieger dilation and Theorem
27 describes functoriality and commutant lifting in this setting.
The proofs are omitted
as they are similar to those appearing in Section 2 of [Po5].

\begin{Proposition} Suppose $\tilde{\uT}$ is the minimal Cuntz-Krieger
dilation of a contractive tuple $\uT$ on a Hilbert space $\CH$ satisfying $A$-relations
with respect to some matrix $A$.
\begin{enumerate}
\item Consider a unital $C^*$-algebra $C^*(\ud)$
generated by the entries of the tuple $\ud=(d_1,\cdots,d_n)$ where
the entries satisfy
$$d^*_id_i=I-\sum^n_{j=1}(1-a_{ij})d_jd^*_j.$$
Assume that $\ud$ also satisfies $ d^*_i d_j = 0$ for $1\leq i \neq j \leq n.$ Let
there be a completely positive map $\varrho: C^*(\ud) \to B(\CH)$ such
that $ \varrho (\ud^\alpha (\ud^\beta)^* )= \uT^\alpha (\uT^\beta)^*.$
Then there is a $*$-homomorphism form $C^*(\ud)$ to $C^*(\tilde{\uT})$
such that $d_i \mapsto \tilde{T}_i$ for all $1 \leq i\leq n.$

\item Suppose $\pi: C^*(\uT) \to B(\tilde{\CK})$ is a
$*$-homomorphism and $\theta : C^* (\tilde{\uT}) \to C^*(\uT)$ the
completely positive map obtained by restricting the compression map (to
$B(\CH)$) of $B(\tilde{\CH})$ to $C^*(\tilde{\uT}).$
Assume the minimal Stinespring dilation of $\pi \circ
\theta$ to be $\tilde{\pi}$ i.e.\:$\pi \circ \theta(X)= P_{\tilde{\CK}} \tilde{\pi}(X) |_{\tilde{\CK}}.$
Then $(\tilde{\pi}(\tilde{T}_1),\cdots,\tilde{\pi}(\tilde{T}_n))$ is the
minimal Cuntz-Krieger dilation of $(\pi(T_1),\cdots,\pi(T_n)).$
\end{enumerate}
\end{Proposition}
\begin{Theorem}
Let $\uT$ be a contractive $n$-tuple on $\CH$ satisfying $A$-relations and $\tilde{\uT}$ be its
minimal Cuntz-Krieger dilation.
\begin{enumerate}
\item Suppose $\pi_1$ and $\pi_2$ are two $*$-homomorphism from
$C^*(\uT)$ to $B(\CK_1)$ and $B(\CK_2)$ respectively, for some Hilbert spaces
$\CK_1$ and $\CK_2.$ Let $\theta$ be as defined in
the previous Proposition. If $X$ is an operator such that $X
\pi_1(Y) = \pi_2(Y) X $ for all $Y \in C^*(\uT),$ and
$\tilde{\pi}_1$ and $\tilde{\pi}_2$ are the minimal Stinespring
dilations of $\pi_1 \circ \theta$ and $\pi_2 \circ \theta$
respectively then there exists an operator $\tilde{X}$ such that
$\tilde{X}\tilde{\pi}_1=\tilde{\pi}_2\tilde{X}$ and
$\tilde{X}P_{\CK_1}=P_{\CK_2}\tilde{X}.$
\item  If $X \in C^*(\uT)'$
then there exists a unique $\tilde{X}\in C^*(\tilde{\uT})' \cap
\{P_\CH\}'$ such that $P_\CH \tilde{X}|_{\CH} =X.$ Also the map $X
\mapsto \tilde{X}$ is a $*$-isomorphism.
\end{enumerate}
\end{Theorem}
Many of the results and arguments for minimal Cuntz-Krieger
dilation in the following part of this section are similar to
those of [DKS] and [DP2] for standard noncommuting dilation.
 Using equation (3.4) we observe that
\begin{eqnarray*}
\tilde{T}^*_j\tilde{T}^*_i \tilde{T}_i
\tilde{T}_j&=&\tilde{T}^*_j[I-\sum_k(1-a_{ik})\tilde{T}_k
\tilde{T}^*_k ] \tilde{T}_j\\
&=&a_{ij} \tilde{T}^*_j\tilde{T}_j \tilde{T}^*_j
\tilde{T}_j=a_{ij} \tilde{T}^*_j\tilde{T}_j.
\end{eqnarray*}
From this follows that for $\alpha=(\alpha_1,\cdots,\alpha_m)$
\begin{equation}
\begin{array}{ccc}
(\tilde{\uT}^\alpha)^* \tilde{\uT}^\alpha &=& a_{\alpha_1
\alpha_2}\cdots
a_{\alpha_{m-1} \alpha_m}\tilde{T}^*_{\alpha_m} \tilde{T}_{\alpha_m}\\
&=& a_{\alpha_1 \alpha_2}\cdots a_{\alpha_{m-1} \alpha_m}
[I-\sum_k(1-a_{\alpha_m k})\tilde{T}_k\tilde{T}^*_k].
\end{array}
\end{equation}
and
$$
\tilde{\uT}^\alpha(\tilde{\uT}^\alpha)^*
\tilde{\uT}^\alpha=\tilde{\uT}^\alpha
$$
so that each
$\tilde{\uT}^\alpha$ is a partial isometry. Let $\hat{\CH}$ and
$\tilde{\CH}$ be the dilation spaces associated with  standard noncommuting dilation $\hat{\uT}$
and $\tilde{\uT}$ respectively as before and let us denote
$\tilde{\CH} \ominus \CH$ by $\CE.$ We know that $\hat{T}_i$ and
$\tilde{T}_i$ leaves $\CE$ and $\hat{\CH} \ominus \CH$  respectively invariant.
Let $\Phi:B(\CH_1) \to B(\CH_1)$ be the completely positive map
defined by
$$ \Phi(X)=\sum^n_{i=1} \hat{T}_i P_{\CH^\bot} X P_{\CH^\bot} \hat{T}^*_i.
$$
Thus, $\Phi (P_\CE) \leq \Phi (I).$ Also let $Q_i:=P_\CE
\tilde{T}_i P_\CE.$ Then for $h \in \CE$
\begin{eqnarray*}
\langle h, \sum_{|\alpha|=m} \uQ^\alpha(\uQ^\alpha)^* h\rangle&= &
\langle  h, \sum_{|\alpha|=m} \tilde{\uT}^\alpha P_{\CE}(\tilde{\uT}
^\alpha)^* h\rangle\\
&=& \langle  h, \sum_{|\alpha|=m} \hat{\uT}^\alpha
P_{\CH^\bot}P_{\CE}P_{\CH^\bot}(\hat{\uT} ^\alpha)^* h
\rangle\\
&=& \langle h,\Phi^m (P_\CE)h \rangle\\
&\leq& \langle h, \Phi^m (I) h \rangle.
\end{eqnarray*}
But as $\lim_{m\to \infty} \langle h,\Phi^m (I) h \rangle=0$ we
have $\lim_{m\to \infty} \langle h,\sum_{|\alpha|=m}
\uQ^\alpha(\uQ^\alpha)^*h\rangle=0$
which implies that $\uQ$ is pure. In the above computation we used
$\hat{T}^*_i$ invariance of $\tilde{\CH}$ for $1\leq i \leq n.$
Here we are interested in understanding the structure of WOT-closed
algebra generated by the minimal Cuntz-Krieger dilation $\tilde{\uT}$ of
some contractive tuple $\uT=(T_1,\cdots,T_n)$ satisfying $A$-relations where $T_i \in
B(\CH).$  Let $\CA$ denote the
WOT-closed algebra generated by all $\tilde{T}_i, 1 \leq i \leq
n$.
\begin{Lemma}
\begin{enumerate}
\item If $\CA$ has no wandering vector then every non-trivial invariant
subspace with respect to $\CA$ is also reducing.
\item $\CK := \CE \ominus[\sum^n_i\tilde{T}_i\CE]$ is a wandering subspace
for $\tilde{\uT}.$
\end{enumerate}
\end{Lemma}
\noindent {\sc Proof:} Let there be no wandering vector for $\CA$
and let if possible $\CN$ be a non-trivial invariant subspace for
$\CA.$ If $\sum^n_{i=1} \tilde{T}_i \CN$  is not equal to $\CN$
then $\CN\ominus \sum^n_{i=1}\tilde{T}_i \CN$ would be wandering
as seen using orthogonality of the ranges of the $\tilde{T}_i$'s, equation
(5.1) and the following: For $n_1, n_2 \in \CN \ominus
\sum^n_{i=1}\tilde{T}_i \CN$
$$ \langle \tilde{T}^*_i \tilde{T}_i \tilde{T}_{\alpha_1}\cdots
\tilde{T}_{\alpha_m}
 n_1, n_2 \rangle= \langle
a_{i \alpha_1}\tilde{T}_{\alpha_1}\cdots \tilde{T}_{\alpha_m} n_1,
n_2 \rangle =0.$$
But this is not possible by our assumption. So
\begin{equation}
\CN= \sum^n_{i=1} \tilde{T}_i \CN.
\end{equation}
Now let $h \in \CN$ be arbitrary. From the above equation it
follows that one can write $h$ as $\sum^n_{i,j=1} \tilde{T}_i
\tilde{T}_j n_{ij}$ for some $n_{ij} \in \CN.$ From this and
equation (5.1) it is clear that $\tilde{T}^*_k h \in \CN$ for all $1
\leq k \leq n.$ So $\CN$ is reducing for $\CA$. Hence (1) follows.

$\CE$ is also an invariant subspace for $\CA.$ The nontrivial case is
when $\CE$ is non-zero. $\CE \neq \sum^n_{i=1} \tilde{T}_i
\CE$ as otherwise $\CE$ would be reducing which is not possible as
$\CH$ spans $\tilde{\CH}.$ It follows from above that $\CK$ is a
wandering subspace of $\CA$. \qed

\medskip

So we can write $\tilde{\CH}=\CH \oplus \CH' \oplus (\Gamma_A \otimes \CK)$
for some Hilbert space $\CH'$. So $\sum_{\alpha \in
\tilde{\Lambda}}\tilde{\uT}^\alpha \CK=\Gamma_A \otimes \CK$
and this is left invariant by all $\tilde{T}_i$. Also
$\tilde{T}_iP_{\Gamma_A \otimes \CK}$ is $ S_i
\otimes I$ for $1 \leq i \leq n.$

Let us denote by $\CB$ the WOT-closed algebra generated by $T_1,
\cdots, T_n.$ In order to get reducing subspaces for $\CA$ it's
sufficient to demand for $\CB^*$-invariant subspace as seen in the
next Lemma. ($\CA[\CL]$ denotes the closed linear span of $A\CL$.)

\begin{Lemma} Let $\CL$ be a $\CB^*$-invariant subspace of $\CH.$
Then $\CA[\CL]$ reduces $\CA.$ If $\CL_1$ and $\CL_2$ are
orthogonal $\CB^*$-invariant subspace of $\CH$ then $\CA[\CL_1]$
and $\CA[\CL_2]$ are also mutually orthogonal.
\end{Lemma}
\noindent {\sc Proof:} $\tilde{T}^*_i$ leaves $\CL$ invariant as
$\tilde{T}_i^*$ and $T^*_i$ leaves $\CH$ and $\CL$ respectively invariant. Thus
$$ \CA[\CL]=\overline{\mbox{span}}\{ \tilde{\uT}^\alpha h : \alpha \in
\tilde{\Lambda}, h \in \CL\}.$$ Now for any $x \in \CL$ and
$\alpha = (\alpha_1, \cdots, \alpha_m),$ using equation (5.1)
\[ \tilde{T}^*_i \tilde{\uT}^\alpha x=\left\{ \begin{array}{cc}
[I-\sum_k (1-a_{ik})\tilde{T}_k\tilde{T}^*_k]x & \mbox{if~}
\alpha_1=i,|\alpha|=1\\
a_{\alpha_1 \alpha_2}\tilde{T}_{\alpha_2} \cdots \tilde{T}_{\alpha_m} x&
\mbox{if~} \alpha_1=i,|\alpha|>1\\
0 & \mbox{if~} \alpha_1 \neq i \\
\tilde{T}^*_ix & \mbox{if~} |\alpha|=0 \end{array} \right. \]

As $\CL$ is invariant for $\CA^*,~~\tilde{T}^*_i\tilde{\uT}^\alpha
x \in \CA[\CL]$ and hence reduce $\CA.$

Further when $\CL_1$ and $\CL_2$ are orthogonal $\CB^*$-invariant
subspaces, to establish that $\CA[\CL_1]$ is $\CA[\CL_2]$ are
orthogonal it is sufficient to check if $|\alpha|\leq |\beta|$
and $l_1 \in \CL_1, l_2 \in \CL_2$, then $\langle
\tilde{\uT}^\alpha l_1, \tilde{\uT}^\beta l_2 \rangle = 0.$ This
is checked easily for all cases by orthogonality of ranges
of different $\tilde{T}_i$'s, $\CB^*$-invariance of $\CL_i$ and
the equation (5.1) except  $\alpha=(\alpha_1,\cdots,\alpha_m)=\beta$. In this case
\begin{eqnarray*}
 \langle (\tilde{\uT}^\alpha)^*(\tilde{\uT}^\alpha)
l_1, l_2 \rangle &=& \langle a_{\alpha_1 \alpha_2}\cdots
a_{\alpha_{m-1} \alpha_m} [I-\sum_k(1-a_{\alpha_m
k})\tilde{T}_k\tilde{T}^*_k] l_1, l_2 \rangle \\
&=&  a_{\alpha_1 \alpha_2}\cdots a_{\alpha_{m-1} \alpha_m}
\{\langle l_1, l_2 \rangle -\sum_k (1-a_{\alpha_m
k})\langle\tilde{T}^*_k l_1,\tilde{T}^*_k l_2 \rangle \}  =0,
\end{eqnarray*}
hence the Lemma
follows. \qed

\medskip

Recall that $\tilde{\CH}_N$ denotes the summand on which the compact operators
in $C^*(\tilde{\uT})$ act trivially. Let  $\CH_N:=\tilde{\CH}_N \cap \CH.$  In the next Proposition we
assume $\CH$ to be finite dimensional.

\begin{Proposition} Let $\uT$ be a contractive tuple satisfying $A$-relations
of operators on a finite dimensional Hilbert space $\CH.$
\begin{enumerate}
\item Let $\CK$ be a reducing subspace of $\tilde{\CH}_N$ with respect to
$\CA$ and let  $h\in \tilde{\CH}$ such that $P_\CK h $
is non-zero. Then there exists $k \in \CA^*[h] \cap \H_N$ such
that $P_\CK k$ is non-zero.
\item Any non-zero subspace of $\tilde{\CH}_N$ which is co-invariant with
respect to
$\tilde{T}_i, 1 \leq i \leq n$ has a non-trivial intersection with
$\CH_N.$
\item $\tilde{\CH}_N=\CA[\CH_N].$ When $\sum T_i T^*_i=I$ and
$\CB=B(\CH)$ every co-invariant subspace of $\CH$ with respect to
all $\tilde{T}_i$'s contains $\CH.$
\end{enumerate}
\end{Proposition}
\noindent {\sc Proof:} Proof is similar to the proof of Lemma 4.1,
Corollary 4.2 and Corollary 4.3 in [DPS]. It uses 
above two Lemmas, pureness of $\uQ,$ Wold decomposition of $\tilde{T}$ and 
compactness of the unit ball of finite dimensional Hilbert space $\CH.$ 
(2) and (3) are corollaries of (1). \qed

\medskip

Now we consider the tuple $\uR$ consisting of right creation operators
on $\Gamma  (\C^n)$ given by $R_i x = x \otimes e_i$.
One can easily notice using methods similar to proof of
Lemma 4 that for polynomials
$p_{(l,m)}=z_lz_m-a_{ml}z_lz_m, (l,m) \in \{1, \cdots, n \}\times
\{1,\cdots, n\} =\CI,$ we get $(\Gamma
(\C^n))^p(\uR)=\Gamma_A.$
 Let $X_i$ denote the compression of $R_i$ to $\Gamma_A$, i. e.
$X_i=P_{\Gamma_A}R_i|_{\Gamma_A}.$
Suppose $e^\alpha \in \Gamma_A,$ and when $ |\alpha| > 0$
let $\alpha=(\alpha_1, \cdots, \alpha_m).$ Then
$$X_i e^\alpha = P_{\Gamma_A} R_ie^\alpha= \left\{
\begin{array}{cc}
e_i & \mbox{~if~} |\alpha|=0 \\
a_{\alpha_m i }   e^\alpha \otimes e_i& \mbox{~if~} |\alpha| \geq
1
\end{array} \right.$$
Moreover from Proposition 9 it follows that $\uX$ consists of isometries with
orthogonal ranges satisfying  $A^t$-relations, where $A^t$ is the transpose of $A$.
Let $\CS$ and $\CX$ denote the WOT-closed algebras generated by
$S_1,\cdots,S_n$ and $X_1,\cdots, X_n$ respectively. Now we shall
analyze the structure of these WOT-closed algebras. Let $Q_k$
denote the projection onto $\mbox{span}\{ e^\alpha:  \alpha\in
\tilde{\Lambda}_A,|\alpha|=k\}.$

\begin{Proposition}
\begin{enumerate}
\item $\CS$ coincides with the commutant of $\CX$ in
$B(\Gamma_A),$ that is $\CS=\CX'.$ Also $\CX=\CS'$ and hence $\CS$
and $\CX$ are double commutants of themselves.
\item $\CS$ and $\CX$ are inverse closed and the only normal
elements in
$\CS$ and $\CX$ are scalars.
\end{enumerate}
\end{Proposition}
\noindent {\sc Proof:} Any element in $\CS$ can be written as
a formal sum $\sum_{\alpha} b_{\alpha}\uS^\alpha$, where $b_{\alpha}\in \C$ are given by $X\omega = \sum b_{\alpha}e^{\alpha}$. Let
for $\beta = (\beta_1, \cdots, \beta_m),$ $\beta'$ denote $(\beta_m,\cdots,\beta_1).$
$$
\uS^\alpha \uX^{\beta'} e^\gamma=\left\{ \begin{array}{cc}
a_{\alpha_{|\alpha|}\gamma_1}a_{\gamma_{|\gamma|}\beta_1}e^\alpha\otimes
e^{\gamma}\otimes e^{\beta} & \mbox{if~}|\gamma|>0\\
a_{\alpha_{|\alpha|}\beta_1}e^\alpha\otimes e^{\beta} & \mbox{if~}|\gamma|=0
\end{array}\right. =\uX^{\beta'}\uS^\alpha e^{\gamma}.
$$
So, $\CS \subseteq \CX'.$ The converse is similar to the proof of Theorem 1.2
in [DP2] after noticing that $X_i Q_k=Q_{k+1}X_i$ and considering
the Ces\'aro sums
$$
p_k(L)=\sum_{|\alpha|<k}\Big( 1-\frac{|\alpha|}{k} \Big) d_\alpha \uS^\alpha
$$
for $L\omega=\sum_{\alpha \in \tilde{\Lambda}_A}d_\alpha e^\alpha.$

$\CS$ and $\CX$ are  inverse closed as this is the case for any
algebra which is a commutant. This and (2) can be proved by taking the
same approach as that of the proof of Corollary 1.4, 1.5 in [DP2].  \qed

\begin{Proposition}
Any element $A \in \CL$ leaves the range of $\uX^\alpha (\uX^\alpha)^*$ invariant.
\end{Proposition}
\noindent {\sc Proof:}
Note that one can argue as we did for $\uS^\alpha$ and show that $\uX^\alpha$
are partial isometries. Further as $\CL = \CX'$
\begin{eqnarray*} \uX^\alpha
(\uX^\alpha)^* A \uX^\alpha (\uX^\alpha)^* &=& \uX^\alpha
(\uX^\alpha)^* \uX^\alpha A (\uX^\alpha)^*\\
&=& \uX^\alpha A (\uX^\alpha)^*= A \uX^\alpha (\uX^\alpha)^*,
\end{eqnarray*}
the proposition follows. \qed

\medskip

In these algebras the wandering subspace description is much simpler than
the general case as can be seen from the next result.

\begin{Proposition}
\begin{enumerate}
\item If $\CN$ is an invariant subspace of $\CL$ then $\CM=\CN \ominus
\sum^n_{i=1} S_i \CN$ is a wandering subspace and $\CL [\CM]=\CN.$
\item A subspace is cyclic and invariant with respect to $\CL$ if
and only if it is the range of some element in $\CX.$
\end{enumerate}
\end{Proposition}
\noindent {\sc Proof:} Follows from the Wold decomposition using methods
similar to the proof of Theorem 2.1 in [DP2]. \qed

\vsp \noindent{\bf Acknowledgements:}
The first author is funded by the Department of Science and Technology (India) under the Swarnajayanthi Fellowship scheme. Second author is supported by Deutscher
Akademischer Austausch Dienst Fellowship.

\end{section}

\begin{center}
{\bf References}
\end{center}
\begin {itemize}

\item [{[AP]}] Arias, A.;  Popescu, G.: Noncommutative interpolation and
Poisson transforms, Israel J. Math., {\bf 115}(2000), 205-234,
{\bf MR 2001i:47021}.

\item [ {[Ar1]}]   Arveson, W.: Subalgebras of $C^*$-algebras. Acta
Math. 123 (1969) 141-224.  {\bf MR 40 \#6274}.

\item [{[Ar2]}] Arveson, W.: {\em An Invitation to $C^*$-algebras,\/ }
Graduate Texts in Mathematics, No. 39, Springer-Verlag, New
York-Heidelberg (1976). {\bf MR 58$\#$23621}.

\item [{[Ar3]}]  Arveson, W.:   Subalgebras of $C^*$-algebras III,
Multivariable
operator theory, Acta Math., {\bf 181}(1998), no. 2, 159-228. {\bf
MR 2000e:47013}.

\item [ { [BB]}]  Bhat, B. V. R.; Bhattacharyya, T.: A model theory for $q$-commuting contractive tuples,  
J. Operator Theory  {\bf 47}  (2002),  no. 1, 97-116. {\bf MR 2003c:47018}.

\item [ { [BBD]}] Bhat, B. V. R.;   Bhattacharyya, T.; Dey,
S. : Standard noncommuting and commuting dilations of
commuting tuples, Trans. Amer. Math. Soc., {\bf 356} (2004), 1551-1568. {\bf
MR 2005b:47011}

\item [{[Bu]}] Bunce, J. W.: Models for $n$-tuples of noncommuting
operators, J. Funct. Anal., {\bf 57} (1984), 21-30. {\bf MR
85k:47019}.

\item [{[Cu]}] Cuntz, J. : Simple $C^*$-algebras generated by isometries,
Commun. Math. Phys., {\bf 57} (1977),173-185. {\bf MR 57$\#
$7189}.

\item [{[CK]}] Cuntz, J.; Krieger, W.  : A class of $C^*$-algebras
and topological Markov chains., Invent. Math, {\bf 63 }
(1981), 25-40. {\bf MR 82f:46073a}.

\item [{[Da]}] Davis, C.: Some dilation and representation theorems.
Proceedings of the Second International Symposium in West Africa
on Functional Analysis and its Applications (Kumasi, 1979),
159-182. {\bf MR 84e:47012}.

\item [ {[De]}] Dey, S.: Standard dilations of $q$-commuting tuples,
Indian Statistical Institute, Bangalore preprint (2003).

\item [{[DKS]}] Davidson K. R.; Kribs, D. W.; Shpigel, M.E.:
Isometric dilations of non-commuting finite rank $n$-tuples,
Canad. J. Math., {\bf 53} (2001) 506-545. {\bf MR 2002f:47010}.

\item [ {[DP1]}]  Davidson, K. R.  and Pitts,  D. R.: The
algebraic structure of non-commutative analytic Toeplitz algebras
Math. Ann., {\bf 311} (1998) 275-303. {\bf MR 2001c:47082}.

\item [ {[DP2]}] Davidson, K. R.  and   Pitts,  D. R.: Invariant
subspaces and hyper-reflexivity for free semi-group algebras,
Proc. London Math. Soc., {\bf 78} (1999) 401-430. {\bf MR 2000k:47005}.

\item [{[Fr]}]  Frazho, A. E.: Models for noncommuting operators, J.
Funct. Anal., {\bf 48} (1982), 1-11. {\bf MR  84h:47010}.

\item [{[Mu]}] Muhly, P. S.: A finite-dimensional introduction to operator
algebra, Operator algebras and applications (Samos, 1996), 313-354,
NATO Adv. Sci. Inst. Ser. C Math. Phys. Sci., {\bf 495,}
Kluwer Acad. Publ., Dordrecht, 1997.  {\bf MR 98h:46062}.

\item [{[Po1]}]  Popescu, G.:  Isometric dilations for infinite
sequences of noncommuting operators, Trans. Amer. Math. Soc., {\bf
316} (1989), 523-536. {\bf MR 90c:47006}.

\item [{[Po2]}]
 Popescu, G.: Characteristic functions for infinite
sequences of noncommuting operators, J. Operator Theory, {\bf 22}
(1989), 51-71. {\bf MR 91m:47012}.

\item [{[Po3]}] Popescu, G.: von Neumann inequality for $(B({\CH})^n)_1$.
Math. Scand. {\bf 68} (1991), no. 2, 292-304.
{\bf MR 92k:47073}.

\item [{[Po4]}] Popescu, G.: Multi-analytic operators on Fock spaces.
Math. Ann. {\bf 303} (1995), no. 1, 31-46. {\bf  MR 96k:47049}.

\item [ {[Po5]}] Popescu, G.: Universal operator algebras associated
to contractive sequences of non-commuting operators,  J. London
Math. Soc. (2), {\bf 58} (1998), no. 2, 469-479.  {\bf MR 99m:47054}.

\item [{[Po6]}] Popescu, G.: Poisson transforms on some $C\sp *$-algebras
generated by isometries. J. Funct. Anal. {\bf 161} (1999), 27-61. {\bf
MR 2000m:46117}.

\end{itemize}
\normalsize

\noindent {\sc B. V. Rajarama Bhat,   Indian Statistical Institute,\\ 
R. V. College Post, Bangalore 560059, India.}\\
{\sl bhat@isibang.ac.in} \\
{\sc Santanu Dey, Institut f\"ur Mathematik und Informatik,\\
Ernst-Moritz-Arndt-Universit\"at, Friedrich-Ludwig-Jahn-Str. 15a, \\
17487 Greifswald, Germany.} \\
{\sl dey@uni-greifswald.de}\\
{\sc Joachim Zacharias, School of Mathematical Sciences,\\
University of Nottingham, Nottingham, NG7 2RD UK.}\\
{\sl joachim.zacharias@nottingham.ac.uk }

\end{document}